\documentclass[12pt,a4paper]{article}

\usepackage[left=2.5cm,right=2.5cm,    top=2cm,bottom=2cm,bindingoffset=0cm]{geometry}

\usepackage[french,british]{babel}

\usepackage[T1]{fontenc}

\usepackage{amsmath}

\usepackage{amsfonts}
\usepackage{amsbsy}
\usepackage{amssymb}

\usepackage{amsthm}

\theoremstyle{plain}

\theoremstyle{definition}
\newtheorem{definition}{Definition}

\newtheorem{lemma}{Lemma}
\newtheorem{theorem}{Theorem}
\newtheorem{corollary}{Corollary}

\theoremstyle{Remark}

\theoremstyle{definition}

\begin{document}

\title{\textbf{Extrapolation Problem for Continuous Time Periodically Correlated Isotropic Random Fields}}

\date{}

\maketitle

\noindent Bulletin of Mathematical Sciences and Applications\\
ISSN: 2278-9634, Vol. 19, pp 1-23\\
doi:10.18052/www.scipress.com/BMSA.19.1\\


\author{\textbf{Iryna Golichenko}$^{1}$, \textbf{Oleksandr Masyutka}$^{2}$, \textbf{Mykhailo Moklyachuk}$^{3,*}$ \\\\
$^{1}${Department of Mathematical Analysis and Probability Theory, \\
National Technical University of Ukraine, Kyiv 03056, Ukraine},\\
$^{2}${Department of Mathematics and Theoretical Radiophysics,\\
Taras Shevchenko National University of Kyiv, Kyiv 01601, Ukraine},\\
$^{3}${Department of Probability Theory, Statistics and Actuarial Mathematics, \\
Taras Shevchenko National University of Kyiv, Kyiv 01601, Ukraine, \\
$^{*}$Corresponding Author: Moklyachuk@gmail.com}\\\\
     }

\vspace{2ex}

\vspace{2ex}

\renewcommand{\abstractname}{Abstract}
\begin{abstract}
The problem of optimal linear estimation of functionals
depending on the unknown values of a random field
$\zeta(t,x)$, which is mean-square continuous periodically correlated
 with respect to  time argument $t\in\mathbb R$ and  isotropic on the unit sphere ${S_n}$ with respect to spatial argument $x\in{S_n}$.
Estimates are based on
observations of the field $\zeta(t,x)+\theta(t,x)$ at points $(t,x):t<0,x\in S_{n}$,
 where $\theta(t,x)$ is an
uncorrelated with $\zeta(t,x)$
random field, which is mean-square continuous periodically correlated
 with respect to time argument $t\in\mathbb R$ and isotropic on the sphere ${S_n}$ with respect to spatial argument $x\in{S_n}$.
Formulas for calculating the mean square errors and the spectral characteristics of the optimal linear
estimate of functionals are derived in the case of spectral certainty where the spectral densities of the fields are exactly known.
Formulas that determine the least
favourable spectral densities and the minimax (robust) spectral
characteristics are proposed in the case where the spectral densities are not exactly known while a class of admissible spectral densities is given.
\end{abstract}

\section{Introduction}

Cosmological Principle (first coined by Einstein): the Universe is,
in the large, homogeneous and isotropic (J.~G.~Bartlett \cite{Bartlett}). Last
decades indicate growing interest to the
spatio-temporal data measured on the surface of a sphere. These data
includes cosmic microwave background (CMB) anisotropies (J.~G.~Bartlett \cite{Bartlett}, W.~Hu and S.~Dodelson \cite{Hu}, N.~Kogo and N.~Komatsu \cite{Kogo}, T.~Okamoto and W.~Hu \cite{Okamoto},
 P.~Adshead and W.~Hu \cite{Adshead}), medical imaging (R.~Kakarala \cite{Kakarala}), global
and land-based temperature data (P.~D.~Jones \cite{Jones}, T.~Subba Rao and G.~Terdik\cite{SubbaRao2006}), gravitational and geomagnetic data, climate model (G.~R.~North and R.~F.~Cahalan \cite{North}). Some basic results and references on the theory of isotropic random fields on a sphere can be found in the
books by M.~I.~Yadrenko \cite{Yadrenko} and A.~M.~Yaglom  \cite{Yaglom:1987a, Yaglom:1987b}. For more recent applications and results see new books by C.~Gaetan and X.~Guyon \cite{Gaetan},
N.~Cressie and C.~K.~Wikle \cite{Cressie}, D.~Marinucci and G.~Peccati \cite{Marinucci} and several
papers covering a number of problems in general for spatial temporal isotropic
observations (T.~Subba Rao and G.~Terdik\cite{SubbaRao2012}, G.~Terdik \cite{Terdik2015}).

Periodically correlated processes and fields are not homogeneous but have numerous properties similar to properties of
stationary processes and homogeneous fields. They describe appropriate models of numerous physical
and man-made processes. A comprehensive list of the existing
references up to the year 2005 on periodically correlated processes
and their applications was proposed by E.~Serpedin, F.~Panduru, I.~Sari and G.~B.~Giannakis\cite{Serpedin}. See also reviews by J.~Antoni \cite{Antoni} and A.~Napolitano \cite{Napolitano}. For more
details see a survey paper by W. A. Gardner \cite{Gardner} and book by H. L. Hurd and A. Miamee \cite{Hurd}. Note, that in the literature periodically
correlated processes are named in multiple different ways such as
cyclostationary, periodically nonstationary or cyclic correlated
processes.

The mean square optimal estimation problems for periodically
correlated with respect to time isotropic on a sphere random fields
are natural generalization of the linear extrapolation,
interpolation and filtering problems for stationary stochastic
processes and homogeneous random fields.
Effective methods of solution of the linear extrapolation, interpolation and filtering problems for stationary stochastic processes and random fields were developed under the condition of certainty where spectral densities of processes
and fields are known exactly (see, for example, selected works of A.~N.~Kolmogorov \cite{Kolmogorov}, survey article by T. Kailath \cite{Kailath}, books by Yu.~A.~Rozanov\cite{Rozanov}, N.~Wiener \cite{Wiener}, A.~M.~Yaglom  \cite{Yaglom:1987a, Yaglom:1987b}, M.~I.~Yadrenko \cite{Yadrenko}, articles by M. P. Moklyachuk and M. I. Yadrenko \cite{MoklyachukYadrenko:1979} - \cite{MoklyachukYadrenko:1980}).

The classical approach to the problems of interpolation, extrapolation and filtering of stochastic processes and random fields is based on the assumption that the spectral densities of processes and fields are known. In practice, however, complete information about the spectral density is impossible in most cases. To overcome this complication one finds parametric or nonparametric estimates of the unknown spectral densities or selects these densities by other reasoning. Then applies the classical estimation method provided that the estimated or selected density is the true one. This procedure can result in a significant increasing of the value of error as K.~S.~Vastola and H.~V.~Poor \cite{Vastola} have demonstrated with the help of some examples. This is a reason to search estimates which are optimal for all densities from a certain class of admissible spectral densities. These estimates are called minimax since they minimize the maximal value of the error of estimates.
Such problems arise when considering problems of automatic control theory, coding and signal processing in radar and sonar, pattern recognition problems of speech signals and images.
A comprehensive survey of results up to the year 1985 in minimax (robust) methods of data processing can be found in the paper by S.~A.~ Kassam and H.~V.~ Poor \cite{KassamPoor}. J. Franke \cite{Franke}, J. Franke and H. V. Poor \cite{Franke_Poor} investigated the minimax extrapolation  and filtering problems for stationary sequences with the help of convex optimization methods. This approach makes it possible to find equations that determine the least favorable spectral densities for different classes of densities.
 The paper by Ulf ~Grenander \cite{Grenander} should be marked as the first one where the minimax
approach to extrapolation problem for the functionals from stationary processes was
developed. For more details see, for example,
survey articles  by M. Moklyachuk \cite{Moklyachuk:2015}, M. Luz and M. Moklyachuk \cite{Luz2016},
books by M.~Moklyachuk \cite{Moklyachuk:2008},  M. Moklyachuk and I. Golichenko  \cite{Moklyachuk:2016}, M.~Moklyachuk and O.~Masytka \cite{Moklyachuk:2012}.
In papers by I. I. Dubovets'ka, O.Yu.~Masyutka and M.P. Moklyachuk \cite{Dubovetska9}, \cite{Dubovetska10}, I.~I.~Golichenko, O. Yu. Masyutka, and M. P. Moklyachuk \cite{Dubovetska11}, \cite{Golichenko2} results of investigation of minimax-robust estimation problems for periodically correlated isotropic random fields are
proposed.

In this article we deal with the problem of mean square optimal linear
estimation of the functional
\[
A\zeta ={\int_{0}^{\infty}}{\int_{S_n}} \,\,a(t,x)\zeta
(t,x)\,m_n(dx)dt
\]
which depends on unknown values of a periodically correlated
(cyclostationary with period $T$) with respect to time isotropic on
the unit sphere ${S_n}$ in Euclidean space ${\mathbb E}^n$ random field
$\zeta(t,x)$, $t\ge 0$, $x\in{S_n}$. Estimates are based on
observations of the field $\zeta(t,x)+\theta(t,x)$ at points
$(t,x)$, $t<0$, $x\in{S_n}$, where $\theta(t,x)$ is an
uncorrelated with $\zeta(t,x)$ periodically correlated with respect
to time isotropic on the sphere ${S_n}$ random field. Formulas are
derived for computing the value of the mean-square error and the
spectral characteristic of the optimal linear estimate of the
functional $A\zeta$ in the case of spectral certainty, where spectral densities of the
fields are known. Formulas are proposed that determine the least
favourable spectral densities and the minimax-robust spectral
characteristic of the optimal estimate of the functional $A\zeta$
for concrete classes of spectral densities under the condition that
spectral densities are not known exactly while classes $D =D_F \times D_G$ of
admissible spectral densities are given.

\begin{section}{Spectral properties of continuous time periodically correlated isotropic on
a sphere random fields}

Let $S_{n}$ be a unit sphere in the $n$-dimensional Euclidean space ${\mathbb E}^n$, let $m_{n}(dx)$ be the Lebesgue measure on  $S_{n}$, and let
\[S_{m}^{l}(x),\, l=1, ... , h(m,n);\, m=0, 1, ...\]
be the orthonormal spherical harmonics of degree $m$, where $h(m,n)$ is the number of orthonormal spherical harmonics (see books by A. Erdelyi et al. \cite{Erdelyi} and C. M\"uller \cite{Muller} for more details).

A mean-square continuous random field ${\zeta(t,x)}$, $t\in\mathbb R$, $x\in{S_n}$, with values ${\zeta(t,x)}$ in $H=L_2(\Omega, \cal F, \mathbb P)$, where $L_2(\Omega, \cal F, \mathbb P)$ denotes the Hilbert space of random variables $\zeta$ with zero first moment, ${\mathbb E}{\zeta}=0$, finite second moment, ${\mathbb E}|{\zeta}|^2<\infty$, and the inner product $\langle \zeta,\theta\rangle= {\mathbb E}{\zeta}\overline{\theta}$, is called periodically correlated (cyclostationary with period $T$) isotropic on the sphere $S_{n}$ if for all $t,s\in\mathbb R$ and $x,y\in{S_n}$ the following property holds true
\[
{\mathbb E}\left({\zeta(t+T,x)} \overline
{\zeta (s+T,y)}\right)=B\left(t,s,\cos\vartheta \right),
\]
where $\cos\vartheta=(x,y)$, $\vartheta$ is the angular distance between points $x,y\in{S_n}$.

The correlation function $B\left(t,s,\cos\vartheta \right)$ of the mean-square continuous random field $\zeta(t,x)$ is continuous. It can be represented in the form of the series
\[
B\left(t,s,\cos\vartheta \right)=\frac{1}{\omega_n}\sum_{m=0}^{\infty}h(m,n) \frac
{C_{m}^{(n-2)/2}(\cos\vartheta)} {C_{m}^{(n-2)/2}(1)}\,\,
B^{\zeta}_m(t,s),
\]
where $\omega_n=(2\pi)^{n/2}\Gamma(n/2)$, $C_m^l(z)$ are the Gegenbauer polynomials (see book by M.~I.~Yadrenko \cite{Yadrenko}).

It follows from the Karhunen theorem that the random field $\zeta(t,x)$ itself can be represented in the form of the mean square convergent series (see  K.~Karhunen \cite{Karhunen}, I. I. Gikhman and A. V. Skorokhod \cite{Gikhman})
\begin{equation} \label{II.1.1}
 {\zeta(t,x)}=
{\sum_{m=0}^{\infty}} {\sum_{l=1}^{h(m,n)}}  S_m^l (x)\zeta_{m}^l
(t),
\end{equation}
where
\[ \zeta_{m}^l (t)= {\int_{S_n}}{\zeta(t,x)} S_m^l
(x)\,m_n(dx).
\]
In this representation
\[\zeta_{m}^l (t),\,\, l=1,\ldots,h(m,n);\,\, t\in\mathbb R, m=0,1,\dots
\]
are mutually uncorrelated periodically correlated stochastic processes with the correlation functions $B^{\zeta}_m(t,s)$:
\[
{\mathbb E}\left(\zeta_{m}^l (t+T)\overline{\zeta_{u}^v
(s+T)}\right)=\delta_m^u \delta_l^v\,\,B^{\zeta}_m(t,s),\]
\[l,v=1,\ldots,h(m,n);\,\, m,u=0,1,\dots;\,\,t,s\in\mathbb R,
\]
where $\delta_l^v$ are the Kronecker delta-functions.

Consider two  mutually uncorrelated periodically correlated isotropic random fields  $\zeta(t,x)$ $t\in\mathbb R$, $x\in{S_n}$, and $\theta(t,x)$, $t\in\mathbb R$, $x\in{S_n}$. We construct the following sequences of stochastic functions
\begin{equation} \label{II.1.2}
\{\zeta_{m}^{l}(j, u)=\zeta_{m}^{l}(u+jT),u\in [0,T), j\in\mathbb Z\},
\end{equation}
\begin{equation} \label{II.1.3}
\{\theta_{m}^{l}(j, u)=\theta_{m}^{l}(u+jT),u\in [0,T), j\in\mathbb Z\}
\end{equation}
which correspond to the random fields $\zeta(t,x)$ and $\theta(t,x)$.
The sequences (\ref{II.1.2}) and (\ref{II.1.3})  form the $L_2([0,T);H)$-valued stationary sequences $\{\zeta_{m}^{l}(j), j\in\mathbb
Z\}$ and $\{\theta_{m}^{l}(j),j\in\mathbb Z\}$, respectively, with the correlation functions
\[
R_{m}^{\zeta}(k,j)=
\int_0^T{\mathbb E} \left[\zeta_{m}^{l}(u+kT)\overline{\zeta_{m}^{l}(u+jT)}\right]du
\]
\[
=\int_0^T B_{m}^{\zeta}(u+(k-j)T,u)du =
 R_{m}^{\zeta}(k-j), \]
 \[
R_{m}^{\theta}(k,j)= \int_0^T{\mathbb E} \left[\theta_{m}^{l}(u+kT)\overline{\theta_{m}^{l}(u+jT)}\right]du
\]
\[
=\int_0^T B_{m}^{\theta}(u+(k-j)T,u)du =
 R_{m}^{\theta}(k-j).\]

To describe properties of the stationary sequences $\{\zeta_{m}^{l}(j), j\in\mathbb
Z\}$ and $\{\theta_{m}^{l}(j)$, $j\in\mathbb Z\}$ we define in the space $L_2([0,T);\mathbb{R})$ the following orthonormal basis
\[
\left\{\widetilde{e}_k=\frac{1}{\sqrt{T}}e^{2\pi
i\{(-1)^k\left[\frac{k}{2}\right]\}u/T}, k=1,2,\dots\right\}, \; \langle \widetilde{e}_j,\widetilde{e}_k\rangle=\delta_k^j.\]

Making use of the introduced basis the stationary sequences $\{\zeta_{m}^{l}(j), j\in\mathbb Z\}$ and $\{\theta_{m}^{l}(j),j\in\mathbb Z\}$ can be represented as follows
\begin{equation} \label{II.1.4}
\zeta_{m}^{l}(j)= \sum_{k=1}^\infty \zeta_{mk}^{l}(j)\widetilde{e}_k,
\end{equation}
\[\zeta_{mk}^{l}(j)=\langle\zeta_{m}^{l}(j),\widetilde{e}_k\rangle = \frac{1}{\sqrt{T}}
\int_0^T \zeta_{m}^{l}(j,v)e^{-2\pi i\{(-1)^k\left[\frac{k}{2}\right]\}v/T}dv,\]
\begin{equation} \label{II.1.5}
\theta_{m}^{l}(j)= \sum_{k=1}^\infty \theta_{mk}^{l}(j)\widetilde{e}_k,
\end{equation}
\[
\theta_{mk}^{l}(j)=\langle\theta_{m}^{l}(j),\widetilde{e}_k\rangle= \frac{1}{\sqrt{T}}
\int_0^T \theta_{m}^{l}(j,v)e^{-2\pi i\{(-1)^k\left[\frac{k}{2}\right]\}v/T}dv.
\]

\noindent
Components  $\{\zeta_{mk}^{l}(j),k=1,2,\dots; j\in\mathbb Z\}$ and $\{\theta_{mk}^{l}(j),k=1,2,\dots;j\in\mathbb Z\}$ of the constructed vector-valued stationary sequences $\{\vec{\zeta}_{m}^{l}(j)=(\zeta_{mk}^{l}(j),k=1,2,\dots)^{\top}, j\in\mathbb Z\}$ and $\{\vec{\theta}_{m}^{l}(j)=(\theta_{mk}^{l}(j),k=1,2,\dots)^{\top},j\in\mathbb Z\}$ have the following properties \cite{Kallianpur},
\cite{Moklyachuk:1981}
\[
\mathbb{E}{\zeta_{mk}^{l}(j)}=0, \quad \|\vec{\zeta}_{m}^{l}(j)\|^2_H=\sum_{k=1}^\infty
\mathbb{E}|\zeta_{mk}^{l}(j)|^2=R_{m}^{\zeta}(0), \]
\[
\mathbb{E}{\zeta_{mk}^{l}(j_1)}\overline{\zeta_{mn}^{l}(j_2)}=\langle
K_m^{\zeta}(j_1-j_2)e_k,e_n\rangle,
\]
\[
\mathbb{E}{\theta_{mk}^{l}(j)}=0, \quad \|\vec{\theta}_{m}^{l}(j)\|^2_H=\sum_{k=1}^\infty
\mathbb{E}|\theta_{mk}^{l}(j)|^2=R_{m}^{\theta}(0), \]
\[
\mathbb{E}{\theta_{mk}^{l}(j_1)}\overline{\theta_{mn}^{l}(j_2)}=\langle
K_m^{\zeta}(j_1-j_2)e_k,e_n\rangle,
\]

\noindent where $\{e_k, k=1,2,\dots\}$ is a basis in the space $\ell_2$.
The correlation functions $K_m^{\zeta}(j)$  and $K_m^{\theta}(j)$ of the stationary sequences $\{\vec{\zeta}_{m}^{l}(j), j\in\mathbb Z\}$ and $\{\vec{\theta}_{m}^{l}(j),j\in\mathbb Z\}$ are correlation operator functions in $\ell_2$.

The vector-valued stationary sequences $\{\vec{\zeta}_{m}^{l}(j)=(\zeta_{mk}^{l}(j),k=1,2,\dots)^{\top}, j\in\mathbb Z\}$ and $\{\vec{\theta}_{m}^{l}(j)=(\theta_{mk}^{l}(j),k=1,2,\dots)^{\top},j\in\mathbb Z\}$ have the spectral density functions
\[F_{m}(\lambda)=\left\{ f_{m}^{kn}(\lambda)\right\}_{k,n = 1 }^{\infty}, \quad G_{m}(\lambda)=\left\{ g_{m}^{kn}(\lambda)\right\}_{k,n = 1 }^{\infty},\]

\noindent which are  operator-valued functions of variable $\lambda\in [-\pi,\pi)$ in the space $\ell_2$
if their correlation functions $K_m^{\zeta}(j)$  and $K_m^{\theta}(j)$ can be represented in the form
\[
\langle
K_m^{\zeta}(j)e_k,e_n\rangle
=\frac{1}{2\pi} \int _{-\pi}^{\pi} e^{ij\lambda}\langle F_{m}(\lambda) {e}_k,{e}_n\rangle d\lambda
=\frac{1}{2\pi} \int _{-\pi}^{\pi} e^{ij\lambda}f_{m}^{kn}(\lambda) d\lambda,
\]
\[
\langle
K_m^{\theta}(j)e_k,e_n\rangle
=\frac{1}{2\pi} \int _{-\pi}^{\pi} e^{ij\lambda}\langle G_{m}(\lambda) {e}_k,{e}_n\rangle d\lambda
=\frac{1}{2\pi} \int _{-\pi}^{\pi} e^{ij\lambda} g_{m}^{kn}(\lambda) d\lambda,
\]

For almost all   $\lambda\in [-\pi,\pi)$ the spectral densities $F_{m}(\lambda)$ and $G_{m}(\lambda)$  are kernel operators with integrable kernel norm
\[
\sum_{k=1}^\infty \frac{1}{2\pi} \int _{-\pi}^\pi \langle F_{m}(\lambda)
e_k, e_k\rangle d\lambda=\sum_{k=1}^\infty\langle K_m^{\zeta}(0)
e_k,e_k\rangle=\|\vec{\zeta}_{m}^{l}(j)\|^2_H=R_{m}^{\zeta}(0),
\]
\[
\sum_{k=1}^\infty \frac{1}{2\pi} \int _{-\pi}^\pi \langle G_{m}(\lambda)
e_k, e_k\rangle d\lambda=\sum_{k=1}^\infty\langle K_m^{\theta}(0)
e_k,e_k\rangle=\|\vec{\theta}_{m}^{l}(j)\|^2_H=R_{m}^{\theta}(0).
\]

In the following sections we explore the described spectral properties of random fields to find solution of the estimation problems.

\end{section}

\begin{section}{Hilbert space projection method of extrapolation of continuous time periodically correlated isotropic random fields}

In this section we deal with the problem of mean square optimal linear
estimation of the functional
\[
A\zeta ={\int_{0}^{\infty}}{\int_{S_n}} \,\,a(t,x)\zeta
(t,x)\,m_n(dx)dt
\]
which depends on the unknown values of a periodically correlated
(cyclostationary with period $T$) isotropic   random field
$\zeta(t,x)$, $t\ge 0$, $x\in{S_n}$. Estimates are based on
observations of the field $\zeta(t,x)+\theta(t,x)$ at points
$(t,x)$, $t<0$, $x\in{S_n}$, where $\theta(t,x)$ is an
uncorrelated with $\zeta(t,x)$ periodically correlated isotropic random field.
Formulas are
derived for computing the value of the mean-square error and the
spectral characteristic of the optimal linear estimate of the
functional $A\zeta$ in the case of spectral certainty, where the spectral densities of the
fields are known. Formulas are proposed that determine the least
favorable spectral densities and the minimax-robust spectral
characteristic of the optimal estimate of the functional $A\zeta$
for concrete classes of spectral densities in the case of spectral uncertainty, where the spectral densities are not known exactly while some classes $D =D_F \times D_G$ of
admissible spectral densities are given.


Consider the problem of the mean square optimal linear estimation of the functional
\[
A\zeta ={\int_{0}^{\infty}}{\int_{S_n}} \,\,a(t,x)\zeta
(t,x)\,m_n(dx)dt
\]
which depends on unknown values of a periodically correlated
isotropic random field
$\zeta(t,x)$, $t\ge 0$, $x\in{S_n}$.
Estimates are based on
observations of the field $\zeta(t,x)+\theta(t,x)$ at points
$(t,x)$, $t<0$, $x\in{S_n}$, where $\theta(t,x)$ is an
uncorrelated with $\zeta(t,x)$ periodically correlated isotropic random field.

It follows from representations \eqref{II.1.1} that the functional $A\zeta$ can be represented in the form
\[A\zeta=\int_{0}^{\infty }\int_{S_{n}}a(t,x)\zeta(t,x)m_{n}(dx)dt=
\]
\[=\sum _{m=0}^{\infty }\sum _{l=1}^{h(m,n)}\int _{0}^{\infty }a_{m}^{l}(t)\zeta_{m}^{l} (t)dt=\]
\[
=\sum _{m=0}^{\infty }\sum _{l=1}^{h(m,n)}\sum_{j=0}^\infty\int_{0}^{T} a_m^l(j,u)\zeta_m^l(j,u)du,
\]
\[ a_{m}^l (t)= {\int_{S_n}}{a(t,x)} S_m^l
(x)\,m_n(dx),
\]
\[
a_m^l(j,u)=a_m^l(u+jT),\,u\in [0,T),
\]
\[
\zeta_m^l(j,u)=\zeta_m^l(u+jT),\, u\in [0,T).
\]

Taking into account the decomposition (\ref{II.1.4}) of stationary sequence $\{\zeta_m^l(j)$, $j\in \mathbb Z\}$, the functional $A\zeta$ can be represented in the following form
\[A\zeta=\sum _{m=0}^{\infty }\sum _{l=1}^{h(m,n)}\sum_{j=0}^\infty\sum_{k=1}^{\infty} a_{mk}^l(j)\zeta_{mk}^l(j)=\]
\[
=\sum _{m=0}^{\infty }\sum _{l=1}^{h(m,n)}\sum_{j=0}^\infty\vec{a}_m^l(j) ^{\top}\vec{\zeta}_m^l(j),
\]
\[
\vec{\zeta}_m^l(j)=(\zeta_{mk}^l(j),k=1,2,\dots)^{\top},\]
\[
\vec{a}_m^l(j)=(a_{mk}^l(j),k=1,2,\dots)^{\top}=\]
\[
=(a_{m1}^l(j),a_{m3}^l(j),a_{m2}^l(j),\dots,a_{m(2k+1)}^l(j),a_{m(2k)}^l(j),\dots)^{\top},
\]
\[
 a_{mk}^l(j)=\langle a_m^l(j),\widetilde{e}_k\rangle =\frac{1}{\sqrt{T}} \int_0^T a_m^l(j,v)e^{-2\pi i\{(-1)^k\left[\frac{k}{2}\right]\}v/T}dv.
\]

We will assume that coefficients $\{\vec a_m^l(j),
j=0,1,\dots\}$ which form this representation satisfy the following conditions
\[
\sum _{m=0}^{\infty }\sum _{l=1}^{h(m,n)}\sum_{j=0}^\infty \|
\vec a_m^l(j)\|<\infty,
\]
\begin{equation} \label{II.2.1}
\sum _{m=0}^{\infty }\sum _{l=1}^{h(m,n)}\sum_{j=0}^\infty (j+1)\|
\vec a_m^l(j)\|^2<\infty,\end{equation}
$$\begin{array} {c}\|\vec a_m^l(j)\|^2=\sum_{k=1}^\infty |a_{mk}^l(j)|^2.\end{array} $$
Under these conditions the functional $A\zeta$ has finite second moment and operators defined below with the help of the coefficients $\{\vec a_m^l(j),
j=0,1,\dots\}$ are compact.

Denote by $L_{2}(F)$ the Hilbert space of complex vector functions
\[h(\lambda)=\left\{ h_{m}^l(\lambda):m=0,1, \ldots; l=1,2,\ldots , h(m,n) \right\},\]
\[
 h_m^l(\lambda)=\left\{ h_{mk}^l \right\}_{k=1}^{\infty},\]
that satisfy condition
\[ \sum _{m=0}^{\infty }\sum _{l=1}^{h(m,n)}\int_{-\pi}^{\pi}(h_m^l(\lambda))^{\top}F_m(\lambda)\overline{h_m^l(\lambda)}d\lambda<\infty.\]

We denote by $L_{2}^{-}(F)$ the subspace of $L_{2}(F)$ generated by the functions
\[e^{ij\lambda}\delta_{k},\delta_{k}=\left\{ \delta_k^n \right\}_{n
= 1 }^{\infty}, k=1,2,\dots,  j<0,\]
where $\delta_k^k=1,\delta_k^n=0$, $k\neq n$.

Every linear estimate  $\hat{A}\zeta$ of the functional $A\zeta$ which is based on observations of the sequence  $\{\vec{\zeta}_m^l(j)+\vec{\theta}_m^l(j), j\in \mathbb Z\}$ at points $j<0$  is defined by the spectral characteristic $h(\lambda)\in L_{2}^{-}(F+G)$ and is of the form
\begin{equation}  \label{II.2.2}
\hat{A}{\zeta}=\sum _{m=0}^{\infty }\sum _{l=1}^{h(m,n)}\int_{-\pi}^{\pi}(h_m^l(\lambda))^{\top}(Z_m^{l\zeta}(d\lambda)+Z_m^{l\theta}(d\lambda)),
\end{equation}

\noindent where $Z_m^{l\zeta}(\Delta)=\{
Z_{mk}^{l\zeta}(\Delta) \}_{k = 1}^{\infty}$ and $Z_m^{l\theta}(\Delta)=\{
Z_{mk}^{l\theta}(\Delta) \}_{k = 1}^{\infty}$ are orthogonal stochastic measures of the sequences  $\{\vec{\zeta}_m^l(j),j\in \mathbb Z\}$ and $\{\vec{\theta}_m^l(j),j\in \mathbb Z\}$ respectively.

Suppose that spectral densities of stationary sequences $\{\vec{\zeta}_m^l(j),j\in \mathbb Z\}, \{\vec{\theta}_m^l(j), j\in \mathbb Z\}$ satisfy the following minimality condition ( see \cite{Rozanov} for more details)
\begin{equation} \label{II.2.3}
\int _{-\pi }^{\pi }\, Tr\, \left[\left(F_{m}(\lambda )+G_{m}(\lambda )\right)^{-1} \right]\, d\lambda  \, <\, \infty .
\end{equation}

The mean square error $\Delta(h;F,G)$ of  the linear estimate $\hat{A}\zeta$  with the spectral characteristic  $h_m^l(\lambda)=\sum_{j=1}^\infty \vec h_m^l(j) e^{-ij\lambda}$
can be represented in the form
\[\Delta(h;F,G)=E|A {\zeta}-\hat{A}
{\zeta}|^{2}=\]
\[=\sum _{m=0}^{\infty }\sum _{l=1}^{h(m,n)}\frac{1}{2\pi} \int_{-\pi}^{\pi} ((A_m^l(\lambda)-h_m^l(\lambda))^{\top}F_m(\lambda){\overline{(A_m^l(\lambda)-h_m^l(\lambda)})}+\]
\[+(h_m^l(\lambda))^{\top}G_m(\lambda)\overline{h_m^l(\lambda)})d\lambda,\; A_m^l(\lambda)=\sum_{j=0}^\infty \vec a_m^l(j) e^{ij\lambda}.\]

The spectral characteristic $h(F,G)$ of the optimal linear estimate $\hat A\zeta$ of the functional minimizes the value of the mean square error

\begin{equation} \label{II.2.4}
\begin{split}
\Delta(F,G)=\Delta(h&(F,G);F,G)=\\
=\mathop {\min }\limits_{h \in
L_{2}^{-}(F+G)} \Delta (h;F,G)&=
\mathop {\min }\limits_{\hat{A}
{\zeta}}E|A{\zeta}-\hat{A} {\zeta}|^{2}.
\end{split}
\end{equation}

With the help of the Hilbert space projection method proposed by A.~N.~Kolmogorov (see, for example, book \cite{Kolmogorov}, p.228-280) we
can find formulas for calculation the mean square error $\Delta(F,G)=\Delta(h(F,G);F,G)$ and
the spectral characteristic $h(F,G)$, which is a solution of the optimization problem (\ref{II.2.4}), of the optimal linear estimate $\hat A\zeta$ of the functional
$A\zeta$. Following the method we find the optimal linear estimate $\hat A\zeta$ as projection of
$A\zeta$ on the closed linear subspace $H^{-}(\zeta+\theta)$ generated  by
values of the field $\zeta(t,x)+\theta(t,x)$ at points
$(t,x)$, $t<0$, $x\in{S_n}$,
in the space $H=L_2(\Omega, \cal F, \mathbb P)$.

This projection is determined by conditions:

1) $\hat A\zeta \in H^{-}(\zeta+\theta)$;

2) $ A\zeta-\hat A\zeta\bot H^{-}(\zeta+\theta)$.

The second condition is satisfied if
for all $m=0,1, \ldots; l=1,2,\ldots , h(m,n)$ and $j=-1,-2,\dots$
\[
\frac{1}{2\pi} \int_{-\pi}^{\pi}
\left[
(A_m^l(\lambda)-h_m^l(\lambda))^{\top}F_m(\lambda)-(h_m^l(\lambda))^{\top}G_m(\lambda)\right]e^{-ij\lambda}d\lambda=0.
\]
These relations mean that for all $m=0,1, \ldots; l=1,2,\ldots , h(m,n)$
\[
(A_m^l(\lambda))^{\top}F_m(\lambda)-(h_m^l(\lambda))^{\top}(F_m(\lambda)+G_m(\lambda))=C_m^l(\lambda),
\]
\[C_m^l(\lambda)=\sum_{j=0}^\infty {c}_m^l(j)e^{ij\lambda},
\]
where $\{{c}_m^l(j), j = 0, 1,\dots\}$ are unknown coefficients. It follows from the indicated
relations that the spectral characteristic $h(F,G)$ of the optimal linear estimate $\hat A\zeta$
is of the form
\[
(h_m^l(F,G))^{\top}= ((A_m^l(\lambda))^{\top}F_m(\lambda) - (C_m^l(\lambda))^{\top})(F_m(\lambda)+G_m(\lambda))^{-1}=
\]
\begin{equation}\label{II.2.5}
 = (A_m^l(\lambda))^{\top}-((A_m^l(\lambda))^{\top}G_m(\lambda)+(C_m^l(\lambda))^{\top}) (F_m(\lambda)+G_m(\lambda))^{-1}.
\end{equation}

The first condition, $\hat A\zeta \in H^{-}(\zeta+\theta)$, is satisfied if for all $m=0,1, \ldots; l=1,2,\ldots , h(m,n)$ and $s=0,1,2,\dots$
\[
\sum _{j=0}^{\infty}\left[
\frac{1}{2\pi}\int_{-\pi}^{\pi}{\left[ { F_m(\lambda)(F_m(\lambda)+G_m(\lambda))^{-1}}
\right]}^{\top}e^{i(j-s)\lambda} d\lambda\right]{a}_m^l(j)=
\]
\begin{equation}\label{II.2.55}
=\sum _{j=0}^{\infty}\left[
\frac{1}{2\pi}\int_{-\pi}^{\pi}{\left[ {(F_m(\lambda)+G_m(\lambda))^{-1}}
\right]}^{\top}e^{i(j-s)\lambda} d\lambda\right]{c}_m^l(j).
\end{equation}

To write these relations in more convenient form we introduce operators $\mathbf{B_m}$, $\mathbf{D_m}$, $\mathbf{R_m}$ determined by matrices
$$
\mathbf{B_m}=\{B_m(j,l)\}_{j,l=0}^\infty,\quad \mathbf{
D_m}=\{D_m(j,l)\}_{j,l=0}^\infty,\quad \mathbf{R_m}=\{R_m(j,l)\}_{j,l=0}^\infty$$
composed with the help of the Fourier coefficients
\[
B_m(j,l)=\frac{1}{2\pi}\int_{-\pi}^{\pi}{\left[ {( F_m(\lambda)+G_m(\lambda))^{-1}}
\right]}^{\top}e^{i(j-l)\lambda} d\lambda,\]
\[D_m(j,l)=\frac{1}{2\pi}\int_{-\pi}^{\pi}{\left[ { F_m(\lambda)(F_m(\lambda)+G_m(\lambda))^{-1}}
\right]}^{\top}e^{i(j-l)\lambda} d\lambda,\]
\[R_m(j,l)=\frac{1}{2\pi}\int_{-\pi}^{\pi}{\left[ { F_m(\lambda)(F_m(\lambda)+G_m(\lambda))^{-1}G_m(\lambda)}
\right]}^{\top}e^{i(j-l)\lambda} d\lambda.
\]
and vectors
\[ {\mathbf{a_m^l}}=\left({a}_m^l(j),j =0,1,2,\dots,\right)^{\top}, \quad {\mathbf{c_m^l}}=\left({c}_m^l(j),j =0,1,2,\dots,\right)^{\top}.
\]
Taking into consideration the introduced operators and vectors we can represent equation
\eqref{II.2.55} in the form
\[
\mathbf{B_m}{\mathbf{c_m^l}}=\mathbf{D_m}{\mathbf{a_m^l}},\quad m=0,1, \ldots; l=1,2,\ldots, h(m,n)
\]
This means that the unknown coefficients ${\mathbf{c_m^l}}=\left({c}_m^l(j),j =0,1,2,\dots\right)^{\top}$ are determined by the equation
\[
{\mathbf{c_m^l}}=\mathbf{B_m^{-1}}\mathbf{D_m}{\mathbf{a_m^l}}.\]

\noindent
It follows from the derived relations that the value of the mean square error
$\Delta(F,G)$ of the optimal linear estimate $\hat A\zeta$ of the functional $ A\zeta$ can
be calculated by the formula
\[\Delta(F,G)=\Delta(h(F,G);F,G)=
\]
\[
=
\sum _{m=0}^{\infty }\sum _{l=1}^{h(m,n)}
\Bigg\{
\frac{1}{2\pi}
\int_{-\pi}^{\pi}
\left[(A_m^l(\lambda))^{\top}G_m(\lambda)+(C_m^l(\lambda))^{\top}\right] (F_m(\lambda)+G_m(\lambda))^{-1}
\times
\]
\[\times
F_m(\lambda)
(F_m(\lambda)+G_m(\lambda))^{-1}
\left[(A_m^l(\lambda))^{\top}G_m(\lambda)+(C_m^l(\lambda))^{\top}\right]^*d\lambda+
\]
\[
+\frac{1}{2\pi}
\int_{-\pi}^{\pi}
\left[(A_m^l(\lambda))^{\top}F_m(\lambda)-(C_m^l(\lambda))^{\top}\right] (F_m(\lambda)+G_m(\lambda))^{-1}
\times
\]
\[\times G_m(\lambda)
(F_m(\lambda)+G_m(\lambda))^{-1}
\left[(A_m^l(\lambda))^{\top}F_m(\lambda)-(C_m^l(\lambda))^{\top}\right]^*d\Bigg\}=
\]
\begin{equation} \label{II.2.6}
=\sum _{m=0}^{\infty }\sum _{l=1}^{h(m,n)}(\langle{{\mathbf{a_m^l}},\mathbf{R_m}{\mathbf{a_m^l}}}\rangle+\langle{{\mathbf{c_m^l}},\mathbf{B_m}{\mathbf{c_m^l}}}\rangle).
\end{equation}

\noindent where $\langle{a,b}\rangle$  is the inner product in the space $\ell_2$.

Let us summarize our results and present them in the form of a theorem.

\begin{theorem} \label{thmII.1}
Let  $\{\zeta(t,x), t\in \mathbb{R}, x\in S_n\}$ and $\{\theta(t,x), t\in \mathbb{R}, x\in S_n\}$ be mutually uncorrelated random fields, which are periodically correlated isotropic on the unit sphere ${S_n}$. Let the stationary sequences  $\{\vec\zeta_m^l(j),j\in\mathbb Z\}$ and $\{\vec\theta_m^l(j),j\in\mathbb Z\}$  constructed with the help of relations (\ref{II.1.2}), (\ref{II.1.3}), respectively, have the spectral densities $F_m(\lambda)$ and $G_m(\lambda)$ that satisfy the minimality condition (\ref{II.2.3}). Let coefficients $\{\vec a_m^l(j), j=0,1,\dots\}$ that determine the functional $A\zeta$ satisfy conditions (\ref{II.2.1}). Then the spectral characteristic $h(F,G)$  and the mean square error $\Delta(F,G)$ of the optimal estimate of the functional $A \zeta$ from observations of the field $\zeta(t,x)+\theta(t,x)$ at points $(t,x)$, $t<0$, $x\in{S_n}$  are given by formulas (\ref{II.2.5}), (\ref{II.2.6}) respectively. The optimal estimate $\hat A\zeta$ of the functional $A\zeta$ is calculated by the formula (\ref{II.2.2}).
\end{theorem}

For the problem of mean square optimal estimation of the functional $A \zeta$ from observations of the field $\zeta(t,x)$ without noise  we have the following corollary.

\begin{corollary} \label{corrII.1}
Let $\{\zeta(t,x), t\in \mathbb{R}, x\in S_n\}$ be a random field, which is periodically correlated  isotropic on the unit sphere ${S_n}$. Let the stationary sequence  $\{\vec\zeta_m^l(j),j\in\mathbb Z\}$ constructed with the help of relations (\ref{II.1.2}) has spectral densities $F_m(\lambda)$ that satisfy the minimality condition
\begin{equation} \label{II.2.7}
\int_{-\pi}^{\pi}{Tr{\left[ {(F_m(\lambda))^{-1}} \right]}}
d\lambda <{\infty}.
\end{equation}
Let coefficients $\{\vec a_m^l(j), j=0,1,\dots\}$ that determine the functional $A\zeta$ satisfy conditions (\ref{II.2.1}). Then the spectral characteristic $h(F)$ and the mean square error $\Delta (F)$ of the optimal linear estimate of the functional $A\zeta$ from observations of the field $\zeta(t,x)$ at points $(t,x)$, $t<0$, $x\in{S_n}$ are determined by formulas
\begin{equation} \label{II.2.8}
(h_m^l(F))^{\top}=(A_m^l(\lambda))^{\top}-(C_m^l(\lambda))^{\top}(F_m(\lambda))^{-1},
\end{equation}
\begin{equation} \label{II.2.9}
\Delta(F)=\sum _{m=0}^{\infty }\sum _{l=1}^{h(m,n)}\langle{{\mathbf{c_m^l}},{\mathbf{a_m^l}}}\rangle,
\end{equation}
where ${\mathbf{c_m^l}}=\left\{\vec {c}_m^l(j) \right\}_{j =0}^\infty=\mathbf{B_m^{-1}}{\mathbf{a_m^l}}$ and matrices $\mathbf{B_m}=\{B_m(j,l)\}_{j,l=0}^\infty$ are defined by elements
\[B_m(j,l)=\frac{1}{2\pi}\int_{-\pi}^{\pi}{\left[ { (F_m(\lambda))^{-1}} \right]}^{\top}e^{i(j-l)\lambda}
d\lambda.\]
\end{corollary}

 Theorem \ref{thmII.1} and Corollary \ref{corrII.1} show us that the Fourier coefficients of some functions from the spectral densities can be used in finding the spectral characteristics and the mean square error of the optimal linear estimates of the functionals of the random fields for problems of extrapolation based on observations of the fields without noise as well as on observations with noise.

To solve the problem of extrapolation of stationary sequences A.~N.~Kolmogorov (see, for example, book \cite{Kolmogorov}, p.272-280) proposed a method based on factorization of the spectral density. This method is suitable for solving the extrapolation problems  based on observations without noise whereas Theorem \ref{thmII.1} describes the method of solving the extrapolation problems  based on observations with a noise.

We apply the indicated  method based on factorization of the spectral density to the problem of estimation of the functional from observations without noise.

Suppose that spectral densities of stationary sequence $\{\vec{\zeta}_m^l(j),j\in\mathbb Z\}$ admit the canonical factorization (see, for example, G.~Kallianpur and V.~Mandrekar \cite{Kallianpur}, M.~P.~Moklyachuk \cite{Moklyachuk:1981})
\begin{equation} \label{II.2.10}
F_m(\lambda)=P_m(\lambda)(P_m(\lambda))^*,\;
P_m(\lambda)=\sum_{u=0}^\infty d_m(u)e^{-iu\lambda},
\end{equation}

\noindent where matrices $d_m(u)=\left\{d_{mk}^r(u)\right\}_{k=\overline{1,\infty}}^{r =\overline{1,M}}$ are coefficients of the canonical representation, $M$ is the multiplicity of $\{\vec{\zeta}_m^l(j),j\in\mathbb Z\}$. In this case the spectral characteristic $h(F)$ and the mean square error $\Delta(F)$ of the optimal estimate $\hat A \zeta$  are determined by formulas

\begin{equation} \label{II.2.11}
h_m^l(F)=A_m^l(\lambda)-(Q_m(\lambda))^{\top}S_m^l(\lambda),
\end{equation}
\begin{equation} \label{II.2.12}
\Delta(F)=\sum _{m=0}^{\infty }\sum _{l=1}^{h(m,n)}\|\mathbf{A_m^ld_m}\|^2,
\end{equation}
where matrices $Q_m(\lambda)$ are defined by equations: $Q_m(\lambda)P_m(\lambda)=I_{M}$,
\[
\|\mathbf{A_m^ld_m}\|^2=\sum_{j=0}^\infty \|(\mathbf{A_m^ld_m})(j)\|^2,\]
\[
 (\mathbf{A_m^ld_m})(j)=\sum_{p=0}^\infty (d_m(p))^{\top}\vec {a}_m^l(p+j),\]
\[S_m^l(\lambda)=\sum_{j=0}^\infty (\mathbf{A_m^ld_m})(j) e^{ij\lambda}.
\]

Let us summarize our results and present them in the form of a theorem.
\begin{theorem}
Let $\{\zeta(t,x), t\in \mathbb{R}, x\in S_n\}$ be a random field, which is periodically correlated isotropic on the unit sphere ${S_n}$. Let the stationary sequence  $\{\vec\zeta_m^l(j),j\in\mathbb Z\}$ constructed with the help of relations (\ref{II.1.2}) have the spectral densities $F_m(\lambda)$ that admit the canonical factorization (\ref{II.2.10}). Let coefficients $\{\vec a_m^l(j), j=0,1,\dots\}$ that determine the functional $A\zeta$ satisfy conditions (\ref{II.2.1}).  Then the spectral characteristic $h(F)$ and the mean square error $\Delta (F)$ of the optimal linear estimate of the functional $A\zeta$ from observations of the field $\zeta(t,x)$ at points $(t,x)$, $t<0$, $x\in{S_n}$ are given by formulas (\ref{II.2.11}), (\ref{II.2.12}).
\end{theorem}


\end{section}

\begin{section}{Minimax-robust method of extrapolation}

Formulas (\ref{II.2.5}) --  (\ref{II.2.12}) for calculating the spectral characteristic and the mean square error of the optimal linear estimate of the functional ${A} \zeta$ can be applied in
the case where the spectral densities $F_m(\lambda)$ and $G_m(\lambda)$ of the stationary sequences $\{\vec{\zeta}_m^l(j),j\in\mathbb Z\}$ and $\{\vec{\theta}_m^l(j),j\in\mathbb Z\}$ constructed with the help of the relations (\ref{II.1.2}), (\ref{II.1.3}), are known.
 If the spectral densities are not  exactly known while a set of admissible densities $D=D_{F}\times D_{G}$ is specified, then the minimax approach to estimation of the functional is reasonable. That is we find the estimate which minimizes the mean square error for all spectral densities from a given set $D=D_{F}\times D_{G}$  simultaneously.

\begin{definition}
For a given class of spectral densities $D=D_{F}\times D_{G}$ the spectral densities $F_m^{0}(\lambda)\in D_{F}$ and $G_m^{0}(\lambda)\in D_{G}$ are called the least favorable in  $D$ for the optimal estimate of functional $A\zeta$ if
\[
\Delta(F^0,G^0)=\Delta(h(F^0,G^0);F^0,G^0)={\max_{\substack{(F,G)\in D}}} \Delta (h(F,G);F,G).
\]
\end{definition}
\begin{definition}
For a given class of spectral densities $D=D_{F}\times D_{G}$ the spectral characteristic $h^0(\lambda)$ of the optimal linear estimate of the functional $A \zeta$ is called minimax-robust if
the following relations hold true
$$\begin{array}{c}
 h^0(\lambda)\in H_{D}=\mathop {\bigcap}\limits_{(F,G)\in D} L_{2}^{-}(F+G),\end{array}$$
$$
 \mathop{\min_{\substack {h\in H_{D}}}} {\max_ {\substack {(F,G)\in D}}}\Delta(h;F,G)
=\mathop {\max_{\substack{(F,G)\in D}}}\Delta(h^0;F,G).
$$
\end{definition}

Taking into account the introduced definitions and relations (\ref{II.2.5}) --  (\ref{II.2.12}) we can verify that the following lemmas hold true.
\begin{lemma}
Spectral densities $F_m^{0}(\lambda)\in D_{F}$ and $G_m^{0}(\lambda)\in D_{G}$ which satisfy the minimality condition (\ref{II.2.3}) are the least favorable in the class $D=D_{F}\times D_{G}$  for the optimal linear estimation of the functional $A\zeta$ from observations of the field $\zeta(t,x)+\theta(t,x)$ at points $(t,x)$, $t<0$, $x\in{S_n}$  if the Fourier coefficients of the functions
\[(F_m^0(\lambda)+G_m^0(\lambda))^{-1},\; F_m^0(\lambda)(F_m^0(\lambda)+G_m^0(\lambda))^{-1},\]
\[F_m^0(\lambda)(F_m^0(\lambda)+G_m^0(\lambda))^{-1}G_m^0(\lambda)\]
determine matrices $\mathbf{B_m^0}, \mathbf{D_m^0}, \mathbf{R_m^0}$ giving a solution of the constrained optimization problem
\begin{equation} \label{II.2.18}
\begin{split}
{\sup_{\substack {(F,G)\in
D}}}\sum _{m=0}^{\infty }\sum _{l=1}^{h(m,n)}\left({\langle\mathbf{a_m^l},\mathbf{R_m}\mathbf{{a_m^l}}}\rangle+\langle{\mathbf{B_m^{-1}}\mathbf{D_m}{\mathbf{a_m^l}},
\mathbf{D_m}{\mathbf{a_m^l}}}\rangle\right)&=\\
=\sum _{m=0}^{\infty }\sum _{l=1}^{h(m,n)}(\langle{{\mathbf{a_m^l}},\mathbf{R_m^0}{\mathbf{a_m^l}}}\rangle+\langle{(\mathbf{B_m^0})^{-1}\mathbf{D_m^0}{\mathbf{a_m^l}},\mathbf{D_m^0}{\mathbf{a_m^l}}}\rangle).&
\end{split}
\end{equation}
\end{lemma}
\begin{lemma}
Spectral densities $F_m^{0}(\lambda)\in D_{F}$ which satisfy the minimality condition (\ref{II.2.7}) are the least favorable in the class  $D_F$ for the optimal linear estimation of the functional $A\zeta$ from observations of the field $\zeta(t,x)$ at points $t<0$,  $x\in S_{n}$ if the Fourier coefficients of the functions $(F_m^0(\lambda))^{-1}$ determine matrices $\mathbf{B_m^0}$ giving  a solution of the constrained optimization problem
\begin{equation} \label{II.2.19}
{\sup_{\substack {F\in
D_F}}}\sum _{m=0}^{\infty }\sum _{l=1}^{h(m,n)}\langle{\mathbf{B_m^{-1}}{\mathbf{a_m^l}},
{\mathbf{a_m^l}}}\rangle=\sum _{m=0}^{\infty }\sum _{l=1}^{h(m,n)}\langle{(\mathbf{B_m^0})^{-1}{\mathbf{a_m^l}},{\mathbf{a_m^l}}}\rangle.
\end{equation}
\end{lemma}

\begin{lemma}
Spectral densities $F_m^{0}(\lambda)\in D_{F}$ which admit the canonical factorization (\ref{II.2.10}) are the least favorable in the class  $D_F$ for the optimal linear estimation of the functional $A\zeta$ from observations of the field $\zeta(t,x)$ at points $t<0$,  $x\in S_{n}$ if the  coefficients of factorizations define a solution of the constrained optimization problem
\begin{equation} \label{II.2.20}
\sum _{m=0}^{\infty }\sum _{l=1}^{h(m,n)}\|\mathbf{A_m^ld_m}\|^2\rightarrow \sup,
\end{equation}
\[F_m(\lambda)=\left(\sum_{u=0}^\infty d_m(u) e^{-iu\lambda} \right)\left(\sum_{u=0}^\infty d_m(u) e^{-iu\lambda} \right)^* \in D_F.\]
\end{lemma}

For more detailed analysis of properties of the least favorable spectral densities and the minimax-robust spectral characteristics we observe that the least favorable spectral densities $F^{0}(\lambda)\in D_{F}$, $G^{0}(\lambda)\in D_{G}$ and the minimax spectral characteristic $h^0=h(F^0,G^0)$ form a saddle point of the function  $\Delta(h;F,G)$ on the set $H_{D}\times D$. The saddle point inequalities
\[
\Delta(h^0;F,G)\leq\Delta(h^0;F^0,G^0)\leq\Delta(h;F^0,G^0), \]
\[
\forall
h\in H_{D}, \quad  \forall F\in D_{F}, \quad  \forall G\in D_{G}
\]
hold if $h^0=h(F^0,G^0)$, $h(F^0,G^0)\in H_{D}$ and $(F^0,G^0)$ is a solution of the constrained optimization problem
\begin{equation} \label{II.2.21}
\Delta(h(F^0,G^0);F,G) \rightarrow {\sup},  \quad (F,G)\in D,
\end{equation}
where the functional
\[
\Delta(h(F^0,G^0);F,G)=
\]
\[=\sum _{m=0}^{\infty }\sum _{l=1}^{h(m,n)}\frac{1}{2\pi
}\int_{-\pi}^{\pi}\left[(A_m^l(\lambda))^{\top}G_m^0(\lambda)+
(C_m^{l0}(\lambda))^{\top}\right] (F_m^0(\lambda)+G_m^0(\lambda))^{-1}\times
\]
\[
\times F_m(\lambda)(F_m^0(\lambda)+G_m^0(\lambda))^{-1}\left[G_m^0(\lambda) \overline{A_m^l(\lambda)}+
\overline{C_m^{l0}(\lambda)}\right]d\lambda+\]
\[+\sum _{m=0}^{\infty }\sum _{l=1}^{h(m,n)}\frac{1}{2\pi
}\int_{-\pi}^{\pi}\left[(A_m^l(\lambda))^{\top}F_m^0(\lambda)-(C_m^{l0}(\lambda))^{\top}\right](F_m^0(\lambda)+G_m^0(\lambda))^{-1}\times
\]
\begin{equation}\label{II.2.22}
\times G_m(\lambda)(F_m^0(\lambda)+G_m^0(\lambda))^{-1}\left[F_m^0(\lambda) \overline{A_m^l(\lambda)}-
\overline{C_m^{l0}(\lambda)}\right]d\lambda.
\end{equation}

The constrained optimization problem (\ref{II.2.21}) is equivalent to the following unconstrained optimization problem
\begin{equation}\label{II.2.23}
\Delta_D(F,G)=-\Delta(h(F^0,G^0);F,G)+\delta((F,G)|D)\rightarrow\inf,
\end{equation}
where $\delta((F,G)|D)$ is the indicator function of the set $D$. Solution $(F^{0}(\lambda), G^0(\lambda))$ to the extremum problem (\ref{II.2.23}) is determined by the condition $0\in\partial \Delta_D(F^{0}, G^{0})$  which is necessary for the point  $(F^{0},G^0)$ to belong to the  set of minimums  of a convex functional. Here
$\partial \Delta_D(F^{0}, G^{0})$ is a subdifferential of the convex functional $\Delta_D(F,G)$ at point $(F,G)=(F^{0}, G^{0})$ (see books by R. T. Rockafellar \cite{Rockafellar}, M.~P.~Moklyachuk \cite{Moklyachuk:2008nonsm}).

The form (\ref{II.2.22}) of the functional $\Delta (h(F^{0} ,G^{0} );F,G)$  is convenient for application the method of Lagrange multipliers for
finding solution to the problem (\ref{II.2.23}).
Making use the method of Lagrange multipliers and the form of
subdifferentials of the indicator functions $\delta((F,G)|D)$
we describe relations that determine the least favorable spectral densities in some special classes
of spectral densities (see books by M.~Moklyachuk \cite{Moklyachuk:2008}, M.~Moklyachuk and O.~Masytka \cite{Moklyachuk:2012}, I. I. Golichenko and M. P. Moklyachuk \cite{Moklyachuk:2016} for more details).

\end{section}

\begin{section}{Least favorable spectral densities in the class $D_\varepsilon\times D_{0}$}

Consider the problem of the minimax estimation of the functional $A\zeta$
depending on the unknown values of the random field
$\{\zeta(t,x), t\in \mathbb R, x\in S_n\}$, which is periodically correlated
 and isotropic on the sphere ${S_n}$.
 Estimates are based on observations of the random field
$\zeta(t,x)+\theta(t,x)$ at points $(t,x):$ $ t<0, x\in{S_n}$.
We consider the problem in the case where the spectral densities
 $F_m(\lambda)$, $G_m(\lambda)$ of stationary sequences $\{\vec \zeta_m^l(j), j\in \mathbb{Z}\}$ and $\{\vec \theta_m^l(j), j\in \mathbb{Z}\}$ which are constructed with the help of relations (\ref{II.1.2}), (\ref{II.1.3}), respectively,
 are not known exactly while there are specified the
following pairs of sets of admissible spectral densities.

The first pair is
\begin{equation*}
\begin{split}
 D_{\varepsilon}^{F1} = \bigg\{F(\lambda )|Tr F_{m}(\lambda )=&(1-\varepsilon)Tr U_{m}(\lambda )+\varepsilon Tr V_{m}(\lambda), \\
&\frac{1}{2\pi\omega_{n} }\sum_{m=0}^{\infty}h(m,n) \int _{-\pi }^{\pi } \text{Tr}\,  F_{m}(\lambda )d\lambda  =p\bigg\},
\end{split}
\end{equation*}
\[ D_{0}^{G1} =\bigg\{G(\lambda )|\frac{1}{2\pi\omega_{n} }\sum_{m=0}^{\infty}h(m,n) \int _{-\pi }^{\pi } \text{Tr}\,  G_{m}(\lambda )d\lambda  =q\bigg\}.\]

The second pair of sets of admissible spectral densities is
\[D_{\varepsilon}^{F2} =\biggl\{F(\lambda )|F_{m}^{kk} (\lambda )=(1-\varepsilon)U_{m}^{kk} (\lambda )+\varepsilon V_{m}^{kk} (\lambda ), \]
\[ \frac{1}{2\pi\omega_{n} }\sum_{m=0}^{\infty}h(m,n) \int _{-\pi }^{\pi }F_{m}^{kk} (\lambda )d\lambda  =p_{k}, k=1,2,\dots\biggr\},\]
\[D_{0}^{G2} =\biggl\{\frac{1}{2\pi\omega_{n} }\sum_{m=0}^{\infty}h(m,n) \int _{-\pi }^{\pi }G_{m}^{kk} (\lambda )d\lambda  =q_{k}, k=1,2,\dots\biggr\}.\]

The third pair of sets of admissible spectral densities is
\begin{equation*}
\begin{split}
D_{\varepsilon}^{F3} =\biggl\{F(\lambda )\vert & \left\langle B_1,F_{m}(\lambda )\right\rangle=(1-\varepsilon)\left\langle B_1 , U_{m}(\lambda )\right\rangle+\varepsilon \left\langle B_1 , V_{m}(\lambda )\right\rangle , \\
&\frac{1}{2\pi\omega_{n} }\sum_{m=0}^{\infty}h(m,n) \int _{-\pi }^{\pi }\left\langle B_{1} , F_{m}(\lambda )\right\rangle d\lambda  =p\biggr\},
\end{split}
\end{equation*}
\[D_{0}^{G3} =\biggl\{\frac{1}{2\pi\omega_{n} }\sum_{m=0}^{\infty}h(m,n) \int _{-\pi }^{\pi }\left\langle B_{2} , G_{m}(\lambda )\right\rangle d\lambda  =q\biggr\}.\]

The forth pair of sets of admissible spectral densities is
\begin{equation*}
\begin{split}
D_{\varepsilon}^{F4} =\biggl\{F(\lambda )|F_{m}(\lambda )=&(1-\varepsilon)U_{m}(\lambda )+\varepsilon V_{m}(\lambda ),\\
&\frac{1}{2\pi \omega_{n}}\sum_{m=0}^{\infty}h(m,n) \int _{-\pi }^{\pi }F_{m}(\lambda )d\lambda  =P\biggr\},
\end{split}
\end{equation*}
\[D_{0}^{G4} =\biggl\{G(\lambda )|\frac{1}{2\pi \omega_{n}}\sum_{m=0}^{\infty}h(m,n) \int _{-\pi }^{\pi }G_{m}(\lambda )d\lambda  =Q\biggr\}.\]

Here $V_{m} (\lambda ), U_{m} (\lambda )$ are given matrices of spectral densities, $p, q, p_{k}, q_{k}, k=1,2,\dots$ are given numbers, $B_{1}, B_{2}, P, Q$ are given positive-definite Hermitian matrices.

From the condition $0\in \partial \Delta _{D} (F^{0} ,G^{0} )$ we find the following equations which determine the least favorable spectral densities for these given sets of admissible spectral densities.

For the first pair $D_{\varepsilon}^{F1}\times D_{0}^{G1}$ we have equations
\begin{equation} \label{II.2.24}
\sum_{l=1}^{h(m,n)}(r_{G}^{0}(\lambda))^{*}r_{G}^{0}(\lambda)=(\alpha_{m} ^{2}+\gamma _{m}(\lambda )) (F_{m}^{0} (\lambda )+G_m^0(\lambda))^{2},
\end{equation}
\begin{equation} \label{II.2.25}
\sum_{l=1}^{h(m,n)}(r_{F}^{0}(\lambda))^{*}r_{F}^{0}(\lambda)=\beta_{m}^{2}(F_{m}^{0} (\lambda )+G_m^0(\lambda))^{2},
\end{equation}
where
\[r_{F}(\lambda)=(A_{m}^{l}(\lambda))^{\top} F_{m}(\lambda )-(C_{m}^{l}(\lambda))^{\top},\]
\[r_{G}(\lambda)= (A_{m}^{l}(\lambda))^{\top} G_{m}(\lambda )+(C_{m}^{l}(\lambda))^{\top},\]
\[\gamma _{m}(\lambda )\le 0\,\, \text{and}\,\, \gamma _{m}(\lambda )=0\,\,\text{if}\,\, \text{Tr}\,  F_{m}^{0} (\lambda )>(1-\varepsilon)\text{Tr}\,  U_{m}(\lambda ),\]
and $\alpha_{m} ^{2} , \beta_{m} ^{2}$ are unknown Lagrange multipliers.

For the second pair $D_{\varepsilon}^{F2}\times D_{0}^{G2}$ we have equations
\begin{equation} \label{II.2.26}
\sum_{l=1}^{h(m,n)}(r_{G}^{0}(\lambda))^{*}r_{G}^{0}(\lambda)=(F_{m}^{0} (\lambda )+G_m^0(\lambda))\left\{(\alpha_{mk}^{2} +\gamma _{mk}(\lambda ))\delta _k^n \right\}_{k,n=1}^{\infty} (F_{m}^{0} (\lambda )+G_m^0(\lambda)),
\end{equation}
\begin{equation} \label{II.2.27}
\sum_{l=1}^{h(m,n)}(r_{F}^{0}(\lambda))^{*}r_{F}^{0}(\lambda)=(F_{m}^{0} (\lambda )+G_m^0(\lambda))\left\{\beta _{mk}^{2}\delta _k^n \right\}_{k,n=1}^{\infty} (F_{m}^{0} (\lambda )+G_m^0(\lambda)),
\end{equation}
where
\[ \gamma _{mk}(\lambda )\le 0\,\,\text{and}\,\, \gamma _{mk}(\lambda )=0\,\,\text{if}\,\, F^{0kk}_{m} (\lambda )>(1-\varepsilon)U^{kk}_{m}(\lambda ),\]
and $\alpha _{mk}^{2} , \beta _{mk}^{2}$ are unknown Lagrange multipliers.

For the third pair $D_{\varepsilon}^{F3}\times D_{0}^{G3}$ we have equations
\begin{equation} \label{II.2.28}
\sum_{l=1}^{h(m,n)}(r_{G}^{0}(\lambda))^{*}r_{G}^{0}(\lambda)=(\alpha_{m}^{2} +\gamma_{m}^{'}(\lambda)) (F_{m}^{0} (\lambda )+G_m^0(\lambda))(B_1)^\top(F_{m}^{0} (\lambda )+G_m^0(\lambda)),
\end{equation}
\begin{equation} \label{II.2.29}
\sum_{l=1}^{h(m,n)}(r_{F}^{0}(\lambda))^{*}r_{F}^{0}(\lambda)=\beta_{m}^{2}(F_{m}^{0} (\lambda )+G_m^0(\lambda))(B_2)^\top(F_{m}^{0} (\lambda )+G_m^0(\lambda)),
\end{equation}
where
\[\gamma_{m}^{'}(\lambda)\le 0\,\,\text{and}\,\, \gamma_{m}^{'}(\lambda)=0\,\,\text{if}\,\, \left\langle B_{1} , F_{m}^{0} (\lambda )\right\rangle >(1-\varepsilon)\left\langle B_{1} , U_{m}(\lambda )\right\rangle ,\]
and $\alpha_{m} ^{2} , \beta_{m} ^{2}$ are unknown Lagrange multipliers.

For the forth pair $D_{\varepsilon}^{F4}\times D_{0}^{G4}$ we have equations
\begin{equation} \label{II.2.30}
\sum_{l=1}^{h(m,n)}(r_{G}^{0}(\lambda))^{*}r_{G}^{0}(\lambda)=(F_{m}^{0} (\lambda )+G_m^0(\lambda)) (\vec{\alpha}_{m}\cdot \vec{\alpha}_{m}^{*}+\Gamma _{m}(\lambda ))(F_{m}^{0} (\lambda )+G_m^0(\lambda)),
\end{equation}
\begin{equation} \label{II.2.31}
\sum_{l=1}^{h(m,n)}(r_{F}^{0}(\lambda))^{*}r_{F}^{0}(\lambda)=(F_{m}^{0} (\lambda )+G_m^0(\lambda))\vec{\beta}_{m}\cdot \vec{\beta}_{m}^{*}(F_{m}^{0} (\lambda )+G_m^0(\lambda)),
\end{equation}
where
$\Gamma _{m}(\lambda )$ are Hermitian matrices,
\[\Gamma _{m}(\lambda )\le 0\,\,\text{and}\,\, \Gamma _{m}(\lambda )=0\,\,\text{if}\,\, F_{m}^{0} (\lambda )>(1-\varepsilon)U_{m}(\lambda ),\]
and $\vec{\alpha}_{m}, \vec{\beta}_{m}$ are unknown Lagrange multipliers.

\begin{theorem}
Let the minimality condition (\ref{II.2.3}) hold true. The least favorable spectral densities  $F_m^0(\lambda)$, $G_m^0(\lambda)$  in the classes $D_0$  for the optimal estimation of the functional $A\zeta$ from observations of the field $\zeta(t,x)+\theta(t,x)$ at points $(t,x)$, $t<0$, $x\in{S_n}$ are determined by relations
\eqref{II.2.24}, \eqref{II.2.25} for the first pair $D_{\varepsilon}^{F1}\times D_{0}^{G1}$ of sets of admissible spectral densities,
\eqref{II.2.26}, \eqref{II.2.27} for the second  pair $D_{\varepsilon}^{F2}\times D_{0}^{G2}$ of sets of admissible spectral densities,
\eqref{II.2.28}, \eqref{II.2.29} for the third  pair $D_{\varepsilon}^{F3}\times D_{0}^{G3}$ of sets of admissible spectral densities,
\eqref{II.2.30}, \eqref{II.2.31} for the fourth pair $D_{\varepsilon}^{F4}\times D_{0}^{G4}$ of sets of admissible spectral densities,
constrained optimization problem  (\ref{II.2.18}) and restrictions  on densities from the corresponding classes $D_{\varepsilon}\times D_0$. The minimax spectral characteristic $h(F^0,G^0)$ of the optimal estimate $\hat A \zeta$ is calculated by (\ref{II.2.5}). The mean square error $\Delta(F^0,G^0)$ is calculated by  (\ref{II.2.6}).
\end{theorem}

For the problem of mean square optimal estimation of the functional $A \zeta$ from observations of the field $\zeta(t,x)$ without noise  we have the following corollary.

\begin{corollary}
Let the minimality condition (\ref{II.2.7}) hold true. The least favorable spectral densities $F_{m}^{0}(\lambda)$ in the classes $D_\varepsilon^{Fk}$, $k=1,2,3,4$, for the optimal linear estimation of the functional  $A\zeta$  from observations of the field $\zeta(t,x)$ at points $t<0$,  $x\in S_{n}$ are determined by the following  equations, respectively,
\begin{equation}
\sum_{l=1}^{h(m,n)}((C_{m}^{l0}(\lambda) )^{\top} )^{*}\cdot(C_{m}^{l0}(\lambda) )^{\top}=(\alpha_{m} ^{2}+\gamma _{m}(\lambda ))(F_{m}^{0} (\lambda ))^{2},
\end{equation}
\begin{equation}
\sum_{l=1}^{h(m,n)}((C_{m}^{l0}(\lambda) )^{\top} )^{*}\cdot(C_{m}^{l0}(\lambda) )^{\top}=F_{m}^{0} (\lambda )\left\{(\alpha _{mk}^{2}+\gamma _{mk}(\lambda ))\delta _k^n \right\}_{k,n=1}^{\infty}F_{m}^{0} (\lambda ),
\end{equation}
\begin{equation}
\sum_{l=1}^{h(m,n)}((C_{m}^{l0}(\lambda) )^{\top} )^{*}\cdot(C_{m}^{l0}(\lambda) )^{\top}=(\alpha_{m} ^{2}+\gamma_{m}^{'}(\lambda)) F_{m}^{0} (\lambda )(B_1)^\top F_{m}^{0} (\lambda ),
\end{equation}
\begin{equation}
\sum_{l=1}^{h(m,n)}((C_{m}^{l0}(\lambda) )^{\top} )^{*}\cdot(C_{m}^{l0}(\lambda) )^{\top}=F_{m}^{0} (\lambda ) (\vec{\alpha}_{m}\cdot \vec{\alpha}_{m}^{*}+\Gamma _{m}(\lambda ))F_{m}^{0} (\lambda ),
\end{equation}
constrained optimization problem  (\ref{II.2.19}) and restrictions  on densities from the corresponding classes  $D_\varepsilon^{Fk}$, $k=1,2,3,4$.
 The minimax spectral characteristic $h(F^0)$ of the optimal estimate $\hat A \zeta$ is calculated by (\ref{II.2.8}). The mean square error $\Delta(F^0)$ is calculated by (\ref{II.2.9}).
\end{corollary}

Consider the problem of optimal linear estimation of the functional $A\zeta$ from observations of the field $\zeta(t,x)$ at points $t<0$, $x\in S_n$ in the case where the spectral densities $F_{m}(\lambda)$ admit the canonical factorization (\ref{II.2.10}).

From the condition  $0\in \partial \Delta _{D} (F^{0})$ we find the following  equations which determine the least favorable spectral densities for the classes $D_{\varepsilon}^{Fk}$, $k=1,2,3,4,$ respectively
\begin{equation} \label{II.2.32}
\sum_{l=1}^{h(m,n)}\left(S_m^{l0}(\lambda)\right)^{\top} \overline{\left(S_m^{l0}(\lambda)\right)}=(\alpha_{m} ^{2}+\gamma _{m} (\lambda ))(P_m^0(\lambda))^{\top} \overline{P_m^0(\lambda)},
\end{equation}
\begin{equation}  \label{II.2.33}
\sum_{l=1}^{h(m,n)}\left(S_m^{l0}(\lambda)\right)^{\top} \overline{\left(S_m^{l0}(\lambda)\right)}=(P_m^0(\lambda))^{\top}\left\{(\alpha _{mk}^{2}+\gamma _{mk}(\lambda )) \delta _k^n \right\}_{k,n=1}^{\infty} \overline{P_m^0(\lambda)},
\end{equation}
\begin{equation}  \label{II.2.34}
\sum_{l=1}^{h(m,n)}\left(S_m^{l0}(\lambda)\right)^{\top} \overline{\left(S_m^{l0}(\lambda)\right)}=(\alpha_{m} ^{2}+\gamma_{m}^{'}(\lambda))(P_m^0(\lambda))^{\top}(B_1)^\top \overline{P_m^0(\lambda)},
\end{equation}
\begin{equation} \label{II.2.35}
\sum_{l=1}^{h(m,n)}\left(S_m^{l0}(\lambda)\right)^{\top} \overline{\left(S_m^{l0}(\lambda)\right)}=(P_m^0(\lambda))^{\top}(\vec{\alpha}_{m}\cdot \vec{\alpha}_{m}^{*}+\Gamma _{m}(\lambda ))\overline{P_m^0(\lambda)}.
\end{equation}
\begin{theorem}
The least favorable spectral densities $F_{m}^{0}(\lambda)$ in the classes $D_{\varepsilon}^{Fk}$, $k=1,2,3,4$, for the optimal linear estimation of the functional  $A\zeta$  from observations of the field $\zeta(t,x)$ for $t<0$,  $x\in S_{n}$ are determined by relations (\ref{II.2.32}) -- (\ref{II.2.35}), respectively, constrained optimization problem  (\ref{II.2.20}) and restrictions  on densities from the corresponding classes $D_{\varepsilon}^{Fk}$, $k=1,2,3,4$. The minimax spectral characteristic of the optimal estimate of the functional $A\zeta$ is calculated by the formula (\ref{II.2.11}). The mean square error $\Delta(F^0)$ is calculated by (\ref{II.2.12}).
\end{theorem}

From the condition  $0\in \partial \Delta _{D} (F^{0})$ we find the following  equations which determine the least favorable spectral densities for the classes $D_0^{Fk}$, $k=1,2,3,4,$ respectively
\begin{equation} \label{II.2.36}
\sum_{l=1}^{h(m,n)}\left(S_m^{l0}(\lambda)\right)^{\top} \overline{\left(S_m^{l0}(\lambda)\right)}=\alpha_m^2(P_m^0(\lambda))^{\top} \overline{P_m^0(\lambda)},
\end{equation}
\begin{equation}  \label{II.2.37}
\sum_{l=1}^{h(m,n)}\left(S_m^{l0}(\lambda)\right)^{\top} \overline{\left(S_m^{l0}(\lambda)\right)}=(P_m^0(\lambda))^{\top}\left\{\alpha _{mk}^{2} \delta _k^n \right\}_{k,n=1}^{\infty} \overline{P_m^0(\lambda)},
\end{equation}
\begin{equation}  \label{II.2.38}
\sum_{l=1}^{h(m,n)}\left(S_m^{l0}(\lambda)\right)^{\top} \overline{\left(S_m^{l0}(\lambda)\right)}=\alpha_m^2(P_m^0(\lambda))^{\top}(B_2)^\top \overline{P_m^0(\lambda)},
\end{equation}
\begin{equation} \label{II.2.39}
\sum_{l=1}^{h(m,n)}\left(S_m^{l0}(\lambda)\right)^{\top} \overline{\left(S_m^{l0}(\lambda)\right)}=(P_m^0(\lambda))^{\top}\vec{\alpha}_{m}\cdot \vec{\alpha}_{m}^{*} \overline{P_m^0(\lambda)}.
\end{equation}
\begin{theorem}
The least favorable spectral densities $F_{m}^{0}(\lambda)$ in the classes $D_{0}^{k}$, $k=1,2,3,4$, for the optimal linear estimation of the functional  $A\zeta$  from observations of the field $\zeta(t,x)$ for $t<0$,  $x\in S_{n}$ are determined by relations (\ref{II.2.36}) -- (\ref{II.2.39}), respectively, constrained optimization problem  (\ref{II.2.20}) and restrictions  on densities from the corresponding classes $D_0^{k}$, $k=1,2,3,4$. The minimax spectral characteristic of the optimal estimate of the functional $A\zeta$ is calculated by the formula (\ref{II.2.11}). The mean square error $\Delta(F^0)$ is calculated by (\ref{II.2.12}).
\end{theorem}

\end{section}

\begin{section}{Least favorable spectral densities in the class $D_V^U\times D_{1\varepsilon}$}

Consider the problem of the minimax estimation of the functional $A\zeta$
depending on the unknown values of the random field
$\{\zeta(t,x), t\in \mathbb R, x\in S_n\}$, which is periodically correlated
 isotropic on the sphere ${S_n}$.
 Estimates are based on observations of the random field
$\zeta(t,x)+\theta(t,x)$ at points $(t,x):$ $ t<0, x\in{S_n}.$
We consider the problem in the case where the spectral densities
 $F_m(\lambda)$, $G_m(\lambda)$ of stationary sequences $\{\vec{\zeta}_m^l(j), j\in \mathbb{Z}\}$ and $\{\vec{\theta}_m^l(j), j\in \mathbb{Z}\}$ which are constructed with the help of relations (\ref{II.1.2}), (\ref{II.1.3}), respectively,
 are not known exactly while there are specified the
following pairs of sets of admissible spectral densities.

The first pair is
\[{D_{V}^{U}} ^{1} =\biggl\{F(\lambda )|\text{Tr}\,  V_{m}(\lambda )\le \text{Tr}\,  (F_{m}(\lambda))\le \text{Tr}\,  (U_{m}(\lambda) ),\]
\[ \frac{1}{2\pi\omega_{n} }\sum_{m=0}^{\infty}h(m,n) \int _{-\pi }^{\pi } \text{Tr}\,  (F_{m}(\lambda) )d\lambda  =p\biggr\},\]
\[ D_{1\varepsilon}^{1} =\bigg\{G(\lambda )| \frac{1}{2\pi\omega_{n} }\sum_{m=0}^{\infty}h(m,n) \int _{-\pi }^{\pi } |\text{Tr}\, (G_{m}(\lambda) )- \text{Tr}\, (G^1_{m}(\lambda))|d\lambda \leq{\varepsilon}\bigg\}.\]

The second pair of sets of admissible spectral densities is
\[{D_{V}^{U}} ^{2} =\biggl\{F(\lambda )|V_{m}^{kk} (\lambda )\le F_{m}^{kk} (\lambda )\le U_{m}^{kk} (\lambda ),\]
\[\frac{1}{2\pi\omega_{n} }\sum_{m=0}^{\infty}h(m,n) \int _{-\pi }^{\pi }F_{m}^{kk} (\lambda )d\lambda  =p_{k} , k=1,2,\dots\biggr\},\]
\[D_{1\varepsilon}^{2} =\biggl\{G(\lambda )|  \frac{1}{2\pi\omega_{n} }\sum_{m=0}^{\infty}h(m,n) \int _{-\pi }^{\pi }
|G_{m}^{kk} (\lambda )-G_{m}^{1kk} (\lambda )|d\lambda \leq{\varepsilon}_{k}, k=1,2,\dots\biggr\}.
\]

The third pair of sets of admissible spectral densities is
\[{D_{V}^{U}}^{3} =\biggl\{F(\lambda )|\left\langle B_1 , V_{m}(\lambda )\right\rangle \le \left\langle B_1 , F_{m}(\lambda )\right\rangle \le \left\langle B_1 , U_{m}(\lambda )\right\rangle,\]
\[ \frac{1}{2\pi\omega_{n} }\sum_{m=0}^{\infty}h(m,n) \int _{-\pi }^{\pi }\left\langle B_1 , F_{m}(\lambda )\right\rangle d\lambda  =p\biggr\},\]
\[D_{1\varepsilon}^{3} =\biggl\{G(\lambda )|
\frac{1}{2\pi\omega_{n} }\sum_{m=0}^{\infty}h(m,n) \int _{-\pi }^{\pi }
|\left\langle B_{2} , G_{m}(\lambda )\right\rangle- \left\langle B_{2} , G^1_{m}(\lambda )\right\rangle|
d\lambda  \leq{\varepsilon}\biggr\}.\]

The forth pair of sets of admissible spectral densities is
\[{D_{V}^{U}} ^{4} =\biggl\{F(\lambda )|V_{m}(\lambda )\le F_{m}(\lambda )\le U_{m}(\lambda ),\, \frac{1}{2\pi \omega_{n}}\sum_{m=0}^{\infty}h(m,n) \int _{-\pi }^{\pi }F_{m}(\lambda )d\lambda  =P\biggr\},\]
\[D_{1\varepsilon}^{4} =\biggl\{G(\lambda )|\frac{1}{2\pi \omega_{n}}\sum_{m=0}^{\infty}h(m,n) \int _{-\pi }^{\pi }
|G_{m}^{kj} (\lambda )-G_{m}^{1kj} (\lambda )|d\lambda \leq{\varepsilon}_{kj}, k,j=1,2,\dots\biggr\}.
\]
Here $V_{m} (\lambda ), U_{m} (\lambda ), G_{m}^1 (\lambda )$ are given matrices of spectral densities, $B_{1}, B_{2}, P$ are given positive-definite Hermitian matrices, $p, p_{k}, k=1,2,\dots$; ${\varepsilon}, {\varepsilon}_{k}, {\varepsilon}_{kj}, k,j=1,2,\dots$ are given numbers.

From the condition $0\in \partial \Delta _{D} (F^{0} ,G^{0} )$ we find the following equations which determine the least favorable spectral densities for these given sets of admissible spectral densities $D_V^U\times D_{1\varepsilon}$.

For the first pair ${D_{V}^{U}} ^{1}\times D_{1\varepsilon}^{1}$ we have equations
\begin{equation} \label{II.2.40}
\sum_{l=1}^{h(m,n)}(r_{G}^{0}(\lambda))^{*}r_{G}^{0}(\lambda)=(\alpha_{m}^{2} +\gamma _{m_{1}} (\lambda )+\gamma _{m_{2}} (\lambda ))(F_{m}^{0} (\lambda )+G_m^0(\lambda))^{2},
\end{equation}
\begin{equation} \label{II.2.41}
\sum_{l=1}^{h(m,n)}(r_{F}^{0}(\lambda))^{*}r_{F}^{0}(\lambda)=
\beta_{m}^{2}\gamma _{m} (\lambda ))(F_{m}^{0} (\lambda )+G_m^0(\lambda))^{2},
\end{equation}

where
\[r_{F}(\lambda)=(A_{m}^{l}(\lambda))^{\top} F_{m}(\lambda )-(C_{m}^{l}(\lambda))^{\top},\]
\[r_{G}(\lambda)= (A_{m}^{l}(\lambda))^{\top} G_{m}(\lambda )+(C_{m}^{l}(\lambda))^{\top},\]
\[\gamma _{m_{1}} (\lambda )\le 0\,\, \text{and}\,\, \gamma _{m_{1}} (\lambda )=0\,\,\text{if}\,\, \text{Tr}\, (F_{m}^{0} (\lambda ))>\text{Tr}\, ( V_{m}(\lambda )),\]
\[\gamma _{m_{2}} (\lambda )\ge 0\,\, \text{and}\,\, \gamma _{m_{2}} (\lambda )=0\,\,\text{if}\,\, \text{Tr}\,  (F_{m}^{0} (\lambda ))<\text{Tr}\,  (U_{m}(\lambda )),\]
\[\gamma _{m} (\lambda )=\text{sign}\,(\text{Tr}\, (G^0_{m}(\lambda) )- \text{Tr}\, (G^1_{m}(\lambda)))
\,\,\text{if}\,],\, (\text{Tr}\, (G^0_{m}(\lambda) )- \text{Tr}\, (G^1_{m}(\lambda)))\not=0,
\]
and $\alpha_{m}^{2}, \beta_{m}^{2}$ are unknown Lagrange multipliers.

For the second pair ${D_{V}^{U}} ^{2}\times D_{1\varepsilon}^{2}$ we have equations
\[\sum_{l=1}^{h(m,n)}(r_{G}^{0}(\lambda))^{*}r_{G}^{0}(\lambda)= (F_{m}^{0} (\lambda )+G_m^0(\lambda))\times\]
\begin{equation}  \label{II.2.42}
\times\left\{(\alpha_{mk}^{2} +\gamma _{m_{1}k} (\lambda )+\gamma _{m_{2}k} (\lambda ))\delta _{k}^n \right\}_{k,n=1}^{\infty} (F_{m}^{0} (\lambda )+G_m^0(\lambda)),
\end{equation}
\begin{equation} \label{II.2.43}
\sum_{l=1}^{h(m,n)}(r_{F}^{0}(\lambda))^{*}r_{F}^{0}(\lambda)=(F_{m}^{0} (\lambda )+G_m^0(\lambda))\left\{(\beta _{mk}^{2} \gamma _{mk} (\lambda ))\delta _k^n \right\}_{k,n=1}^{\infty} (F_{m}^{0} (\lambda )+G_m^0(\lambda)),
\end{equation}
where
\[ \gamma _{m_{1}k} (\lambda )\le 0\,\,\text{and}\,\, \gamma _{m_{1}k} (\lambda )=0\,\,\text{if}\,\, F^{0kk}_{m} (\lambda )>V^{kk}_{m}(\lambda ),\]
\[ \gamma _{m_{2}k} (\lambda )\ge 0\,\,\text{and}\,\,\gamma _{m_{2}k} (\lambda )=0\,\,\text{if}\,\, F^{0kk}_{m} (\lambda )<U^{kk}_{m}(\lambda ),\]
\[\gamma _{mk} (\lambda )=\text{sign}\,(G^{0kk}_{m} - G^{1kk}_{m})
\,\,\text{if}\,\,\,(G^{0kk}_{m} - G^{1kk}_{m})\not=0,
\]
and $\alpha _{mk}^{2} , \beta _{mk}^{2}$ are unknown Lagrange multipliers.

For the third pair ${D_{V}^{U}} ^{3}\times D_{1\varepsilon}^{3}$ we have equations
\begin{equation} \label{II.2.44}
\sum_{l=1}^{h(m,n)}(r_{G}^{0}(\lambda))^{*}r_{G}^{0}(\lambda)=(\alpha_{m}^{2} +\gamma _{m_{1}}^{'}(\lambda )+\gamma _{m_{2}}^{'}(\lambda ))(F_{m}^{0} (\lambda )+G_m^0(\lambda))(B_1)^\top(F_{m}^{0} (\lambda )+G_m^0(\lambda)),
\end{equation}
\begin{equation} \label{II.2.45}
\sum_{l=1}^{h(m,n)}(r_{F}^{0}(\lambda))^{*}r_{F}^{0}(\lambda)=(\beta_{m}^{2} \gamma_{m}^{'}(\lambda)) (F_{m}^{0} (\lambda )+G_m^0(\lambda))(B_2)^\top(F_{m}^{0} (\lambda )+G_m^0(\lambda)),
\end{equation}
where
\[\gamma _{m_{1}}^{'} (\lambda )\le 0\,\,\text{and}\,\, \gamma _{m_{1}}^{'} (\lambda )=0\,\,\text{if}\,\, \left\langle B_{1} , F_{m}^{0} (\lambda )\right\rangle >\left\langle B_{1} , V_{m} (\lambda )\right\rangle ,\]
\[ \gamma _{m_{2}}^{'} (\lambda )\ge 0\,\,\text{and}\,\,\gamma _{m_{2}}^{'} (\lambda )=0 \,\,\text{if}\,\, \left\langle B_{1} , F_{m}^{0} (\lambda )\right\rangle <\left\langle B_{1} , U_{m} (\lambda )\right\rangle ,\]
\[\gamma_{m}^{'}(\lambda)=\text{sign}\,\left(\left\langle B_{2} , G_{m}^{0} (\lambda )\right\rangle -\left\langle B_{2} , G_{m}^{1} (\lambda )\right\rangle\right)
\,\,\text{if}
\]
\[
\left\langle B_{2} , G_{m}^{0} (\lambda )\right\rangle -\left\langle B_{2} , G_{m}^{1} (\lambda )\right\rangle
  \not=0,
\]
and $\alpha_{m} ^{2} , \beta_{m} ^{2}$ are unknown Lagrange multipliers.

For the forth pair ${D_{V}^{U}} ^{4}\times D_{\varepsilon}^{4}$ we have equations
\begin{equation} \label{II.2.46}
\sum_{l=1}^{h(m,n)}(r_{G}^{0}(\lambda))^{*}r_{G}^{0}(\lambda)= (F_{m}^{0} (\lambda )+G_m^0(\lambda))(\vec{\alpha}_m\cdot \vec{\alpha}_m^{*}+\Gamma _{m_{1}} (\lambda )+\Gamma _{m_{2}} (\lambda )) (F_{m}^{0} (\lambda )+G_m^0(\lambda)),
\end{equation}
\begin{equation} \label{II.2.47}
\sum_{l=1}^{h(m,n)}(r_{F}^{0}(\lambda))^{*}r_{F}^{0}(\lambda)=(F_{m}^{0} (\lambda )+G_m^0(\lambda)) (\vec{\beta}_{m}\Gamma _{m} (\lambda ) \vec{\beta}_{m}^{*})(F_{m}^{0} (\lambda )+G_m^0(\lambda)),
\end{equation}
where
$\Gamma _{m_{1}} (\lambda ), \Gamma _{m_{2}} (\lambda ), \Gamma_{m} (\lambda )$ are Hermitian matrices,
\[\Gamma _{m_{1}} (\lambda )\le 0\,\,\text{and}\,\, \Gamma _{m_{1}} (\lambda )=0\,\,\text{if}\,\, F_{m}^{0} (\lambda )>V_{m}(\lambda ),\]
\[\Gamma _{m_{2}} (\lambda )\ge 0\,\,\text{and}\,\, \Gamma _{m_{2}} (\lambda )=0 \,\,\text{if}\,\, F_{m}^{0} (\lambda )<U_{m} (\lambda ),\]
\[\Gamma _{m}^{kj} (\lambda )=
\text{sign}\,(G^{0kj}_{m} - G^{1kj}_{m})
\,\,\text{if}\,\,(G^{0kj}_{m} - G^{1kj}_{m})\not=0,
\]
\noindent and $   \vec{\alpha}_{m}, \vec{\beta}_{m}$ are vectors of unknown Lagrange multipliers.

The following theorem holds true.

\begin{theorem}
Let the minimality condition (\ref{II.2.3})  hold true. The least favorable spectral densities  $F_m^0(\lambda)$, $G_m^0(\lambda)$  in the classes $D_V^U\times D_{1\varepsilon}$  for the optimal estimation of the functional $A\zeta$ from observations of the field $\zeta(t,x)+\theta(t,x)$ at points $(t,x)$, $t<0$, $x\in{S_n}$ are determined by relations
\eqref{II.2.40}, \eqref{II.2.41} for the first pair ${D_{V}^{U}}^{1}\times D_{1\varepsilon}^{1}$ of sets of admissible spectral densities,
\eqref{II.2.42}, \eqref{II.2.43} for the second  pair ${D_{V}^{U}}^{2}\times D_{1\varepsilon}^{2}$ of sets of admissible spectral densities,
\eqref{II.2.44}, \eqref{II.2.45} for the third  pair ${D_{V}^{U}}^{3}\times D_{1\varepsilon}^{3}$ of sets of admissible spectral densities,
\eqref{II.2.46}, \eqref{II.2.47} for the fourth pair ${D_{V}^{U}}^{4}\times D_{1\varepsilon}^{4}$ of sets of admissible spectral densities,
constrained optimization problem   (\ref{II.2.18}) and restrictions  on densities from the corresponding classes $D_V^U\times D_{1\varepsilon}$. The minimax spectral characteristic $h(F^0,G^0)$ of the optimal estimate $\hat A \zeta$ is calculated (\ref{II.2.5}). The mean square error $\Delta(F^0,G^0)$ is calculated by  (\ref{II.2.6}).
\end{theorem}

For the problem of mean square optimal estimation of the functional $A \zeta$ from observations of the field $\zeta(t,x)$ without noise  we have the following corollary.

\begin{corollary}
Let the minimality condition (\ref{II.2.7}) hold true. The least favorable spectral densities $F_{m}^{0}(\lambda)$ in the classes ${D_{V}^{U}} ^{k}$, $k=1,2,3,4$, for the optimal linear estimation of the functional  $A\zeta$  from observations of the field $\zeta(t,x)$ at points $t<0$,  $x\in S_{n}$ are determined by the following  equations, respectively,
\begin{equation}
\sum_{l=1}^{h(m,n)}((C_{m}^{l0}(\lambda) )^{\top} )^{*}\cdot(C_{m}^{l0}(\lambda) )^{\top}=(\alpha_{m}^{2} +\gamma _{m_{1}}^1 (\lambda )+\gamma _{m_{2}}^1 (\lambda )) (F_{m}^{0} (\lambda ))^{2},
\end{equation}
\begin{equation}
\sum_{l=1}^{h(m,n)}((C_{m}^{l0}(\lambda) )^{\top} )^{*}\cdot(C_{m}^{l0}(\lambda) )^{\top}=F_{m}^{0} (\lambda )\left\{(\alpha_{mk}^{2} +\gamma _{m_{1}k}^1 (\lambda )+\gamma _{m_{2}k}^1 (\lambda ))\delta _{k}^n \right\}_{k,n=1}^{\infty}F_{m}^{0} (\lambda ),
\end{equation}
\begin{equation}
\sum_{l=1}^{h(m,n)}((C_{m}^{l0}(\lambda) )^{\top} )^{*}\cdot(C_{m}^{l0}(\lambda) )^{\top}=(\alpha_{m}^{2} +\gamma _{m_{1}}^{1'}(\lambda )+\gamma _{m_{2}}^{1'}(\lambda )) F_{m}^{0} (\lambda )(B_1)^\top F_{m}^{0} (\lambda ),
\end{equation}
\begin{equation}
\sum_{l=1}^{h(m,n)}((C_{m}^{l0}(\lambda) )^{\top} )^{*}\cdot(C_{m}^{l0}(\lambda) )^{\top}=F_{m}^{0} (\lambda )(\vec{\alpha}_m\cdot \vec{\alpha}_m^{*}+\Gamma _{m_{1}}^1 (\lambda )+\Gamma _{m_{2}}^1 (\lambda ))F_{m}^{0} (\lambda ),
\end{equation}
constrained optimization problem  (\ref{II.2.19}) and restrictions  on densities from the corresponding classes ${D_{V}^{U}} ^{k}$, $k=1,2,3,4$.
 The minimax spectral characteristic $h(F^0)$ of the optimal estimate $\hat A\zeta$ is calculated by (\ref{II.2.8}). The mean square error $\Delta(F^0)$ is calculated by (\ref{II.2.9}).
\end{corollary}

\begin{corollary}
Let the minimality condition (\ref{II.2.7}) hold true. The least favorable spectral densities $F_{m}^{0}(\lambda)$ in the classes $D_{1\varepsilon}^{k}$, $k=1,2,3,4$, for the optimal linear estimation of the functional  $A\zeta$  from observations of the field $\zeta(t,x)$ at points $t<0$,  $x\in S_{n}$ are determined by the following  equations, respectively,
\begin{equation}
\sum_{l=1}^{h(m,n)}((C_{m}^{l0}(\lambda) )^{\top} )^{*}\cdot(C_{m}^{l0}(\lambda) )^{\top}=
(\alpha_{m} ^{2}\gamma _{m} (\lambda ))(F_{m}^{0} (\lambda ))^{2},
\end{equation}
\begin{equation}
\sum_{l=1}^{h(m,n)}((C_{m}^{l0}(\lambda) )^{\top} )^{*}\cdot(C_{m}^{l0}(\lambda) )^{\top}=F_{m}^{0} (\lambda )\left\{(\alpha _{mk}^{2}\gamma _{mk} (\lambda ))\delta _k^n \right\}_{k,n=1}^{\infty}F_{m}^{0} (\lambda ),
\end{equation}
\begin{equation}
\sum_{l=1}^{h(m,n)}((C_{m}^{l0}(\lambda) )^{\top} )^{*}\cdot(C_{m}^{l0}(\lambda) )^{\top}=(\alpha_{m} ^{2}\gamma_{m}^{'}(\lambda)) F_{m}^{0} (\lambda )(B_2)^\top F_{m}^{0} (\lambda ),
\end{equation}
\begin{equation}
\sum_{l=1}^{h(m,n)}((C_{m}^{l0}(\lambda) )^{\top} )^{*}\cdot(C_{m}^{l0}(\lambda) )^{\top}=F_{m}^{0} (\lambda ) (\vec{\beta}_{m}\Gamma _{m} (\lambda) \vec{\beta}_{m}^{*})F_{m}^{0} (\lambda ),
\end{equation}

\noindent constrained optimization problem  (\ref{II.2.19}) and restrictions  on densities from the corresponding classes  $D_{1\varepsilon}^{k}$, $k=1,2,3,4$.
 The minimax spectral characteristic $h(F^0)$ of the optimal estimate $\hat A \zeta$ is calculated by formula (\ref{II.2.8}). The mean square error $\Delta(F^0)$ is calculated by (\ref{II.2.9}).
\end{corollary}

Consider the problem of optimal linear estimation of the functional $A\zeta$ from observations of the field $\zeta(t,x)$ at points $t<0$, $x\in S_n$ in the case where the spectral densities $F_{m}(\lambda)$ admit the canonical factorization (\ref{II.2.10}).

From the condition  $0\in \partial \Delta _{D} (F^{0})$ we find the following  equations which determine the least favorable spectral densities for the classes ${D_{V}^{U}} ^{k}, k=1,2,3,4$ respectively
\begin{equation} \label{II.2.48}
\sum_{l=1}^{h(m,n)}\left(S_m^{l0}(\lambda)\right)^{\top} \overline{\left(S_m^{l0}(\lambda)\right)}=(\alpha_{m}^{2} +\gamma _{m_{1}}^1 (\lambda )+\gamma _{m_{2}}^1 (\lambda ))(P_m^0(\lambda))^{\top} \overline{P_m^0(\lambda)},
\end{equation}
\begin{equation}  \label{II.2.49}
\sum_{l=1}^{h(m,n)}\left(S_m^{l0}(\lambda)\right)^{\top} \overline{\left(S_m^{l0}(\lambda)\right)}=(P_m^0(\lambda))^{\top}\left\{(\alpha _{mk}^{2}+\gamma _{m_{1}k}^1 (\lambda )+\gamma _{m_{2}k}^1 (\lambda )) \delta _k^n \right\}_{k,n=1}^{\infty} \overline{P_m^0(\lambda)},
\end{equation}
\begin{equation}  \label{II.2.49}
\sum_{l=1}^{h(m,n)}\left(S_m^{l0}(\lambda)\right)^{\top} \overline{\left(S_m^{l0}(\lambda)\right)}=(\alpha_{m} ^{2}+\gamma _{m_{1}}^{1'}(\lambda )+\gamma _{m_{2}}^{1'}(\lambda ))(P_m^0(\lambda))^{\top}(B_1)^\top \overline{P_m^0(\lambda)},
\end{equation}
\begin{equation} \label{II.2.50}
\sum_{l=1}^{h(m,n)}\left(S_m^{l0}(\lambda)\right)^{\top} \overline{\left(S_m^{l0}(\lambda)\right)}=(P_m^0(\lambda))^{\top}(\vec{\alpha}_{m}\cdot \vec{\alpha}_{m}^{*}+\Gamma _{m_{1}}^1 (\lambda )+\Gamma _{m_{2}}^1 (\lambda ))\overline{P_m^0(\lambda)}.
\end{equation}

\begin{theorem}
The least favorable spectral densities $F_{m}^{0}(\lambda)$ in the classes ${D_{V}^{U}} ^{k}$, $k=1,2,3,4$, for the optimal linear estimation of the functional  $A\zeta$  from observations of the field $\zeta(t,x)$ for $t<0$,  $x\in S_{n}$ are determined by relations (\ref{II.2.48}) -- (\ref{II.2.50}), respectively, constrained optimization problem  (\ref{II.2.20}) and restrictions  on densities from the corresponding classes ${D_{V}^{U}} ^{k}$, $k=1,2,3,4$. The minimax spectral characteristic of the optimal estimate of the functional $A\zeta$ is calculated by the formula (\ref{II.2.11}). The mean square error $\Delta(F^0)$ is calculated by (\ref{II.2.12}).
\end{theorem}

\end{section}

\section{Conclusions}
In this paper we propose formulas for calculating the mean square error and the spectral characteristic of the optimal linear estimate of the
functional
\[
A\zeta ={\int_{0}^{\infty}}{\int_{S_n}} \,\,a(t,x)\zeta
(t,x)\,m_n(dx)dt
\]
depending on unknown values of a mean-square continuous periodically correlated
(cyclostationary with period $T$) with respect to time argument and isotropic on
the unit sphere ${S_n}$ in Euclidean space ${\mathbb E}^n$ random field
$\zeta(t,x)$, $t\in\mathbb R$, $x\in{S_n}$. Estimates are based on
observations of the field $\zeta(t,x)+\theta(t,x)$ at points
$(t,x)$, $t<0$, $x\in{S_n}$ where $\theta(t,x)$ is an
uncorrelated with $\zeta(t,x)$ mean-square continuous periodically correlated with respect
to time argument and isotropic on the sphere ${S_n}$ random field.
The problem is investigated in the case of spectral certainty where matrices of spectral densities of random fields are known exactly and in the case of spectral uncertainty where matrices of spectral densities of random fields are not known exactly while some classes of admissible spectral density matrices are given.
 We derive formulas for calculation the spectral characteristic and the mean-square error of the optimal linear estimate of the functional $A\zeta$
in the case of spectral certainty, where spectral densities $F_{m}(\lambda ), G_{m}(\lambda )$ of the stationary sequences that generate the random fields $\zeta(t,x)$, $\theta (t,x)$ are known exactly.

 We propose a representation of the mean square
error in the form of a linear functional in the $L_1\times L_1$ space with
respect to spectral densities $(F,G)$, which allows us to solve the
corresponding constrained optimization problem and describe the minimax
(robust) estimates of the functional $A\zeta$
for concrete classes of spectral densities under the condition that
spectral densities are not known exactly while classes $D =D_F \times D_G$ of
admissible spectral densities are given.


\begin{thebibliography}{99}


\bibitem{Adshead}
\newblock P. Adshead, and W. Hu,
\newblock \emph{Fast computation of first-order feature-bispectrum corrections},
\newblock Phys. Rev., vol. D85, 103531, 2012.

\bibitem{Antoni}
\newblock J. Antoni,
\newblock \emph{Cyclostationarity by examples},
\newblock Mechanical Systems and Signal Processing, vol. 23, pp. 987--1036, 2009.

\bibitem{Bartlett}
\newblock J. G. Bartlett,
\newblock \emph{The standard cosmological model and cmb anisotropies},
\newblock New Astron. Rev., vol. 43, pp. 83--109, 1999.

\bibitem{Cressie}
\newblock N. Cressie, and C. K. Wikle,
\newblock \emph{Statistics for spatio-temporal data},
\newblock Wiley Series in Probability and Statistics, 2011.




\bibitem{Dubovetska9}
\newblock I. I. Dubovets'ka, O. Yu. Masyutka, and M. P. Moklyachuk,
\newblock \emph{Filtering problems for periodically correlated isotropic random fields},
\newblock Mathematics and Statistics, vol.2, no. 4, pp. 162--171, 2014.

\bibitem{Dubovetska10}
\newblock I. I. Dubovets'ka, O. Yu. Masyutka, and M. P. Moklyachuk,
\newblock \emph{Estimation problems for periodically correlated isotropic random fields},
\newblock Methodology and Computing in Applied Probability, vol.17, no. 1, pp. 41--57, 2015.



\bibitem{Erdelyi}
\newblock A. Erdelyi, W. Magnus, F. Oberhettinger, F. G. Tricomi,
\newblock \emph{Higher transcendental functions. Vol. II.}
\newblock Bateman Manuscript Project. New York-Toronto-London: McGraw-Hill Book Co., Inc. XVII, 1953.



\bibitem{Franke}
\newblock J. Franke,
\newblock \emph{Minimax robust prediction of discrete time series},
\newblock Z. Wahrscheinlichkeitstheor. Verw. Gebiete, vol. 68, pp. 337--364, 1985.

\bibitem{Franke_Poor}
\newblock J. Franke and H. V. Poor,
\newblock \emph{Minimax-robust filtering and finite-length robust predictors},
\newblock Robust and Nonlinear Time Series Analysis. Lecture Notes in Statistics, Springer-Verlag,
vol. 26, pp. 87--126, 1984.


\bibitem{Gaetan}
\newblock C. Gaetan, and X. Guyon,
\newblock \emph{Spatial statistics and modeling},
\newblock Springer Series in Statistics, vol. 81, Springer Science+Business Media, 2010.

\bibitem{Gikhman}
 \newblock I. I. Gikhman and A. V. Skorokhod,
\newblock \emph{The theory of stochastic processes. I.},
\newblock Berlin: Springer, 2004.

\bibitem{Gardner}
\newblock W. A. Gardner,
\newblock \emph{Cyclostationarity in communications and signal processing},
\newblock  New York: IEEE Press, 1994.

\bibitem{Gladyshev}
\newblock E. G. Gladyshev,
\newblock \emph{Periodically correlated random sequences},
\newblock Sov. Math. Dokl. vol. 2, pp. 385--388, 1961.


\bibitem{Dubovetska11}
\newblock I. I. Golichenko, O. Yu. Masyutka, and M. P. Moklyachuk,
\newblock \emph{Minimax-robust fitering of functionals from periodically correlated random fields},
\newblock Cogent Mathematics, vol.2, 1074327, 2015.

\bibitem{Golichenko2}
\newblock I. I. Golichenko, O. Yu. Masyutka, and M. P. Moklyachuk,
\newblock \emph{Filtering of continuous time periodically correlated isotropic random fields},
\newblock Stochastic Modeling and Applications, vol.20, No. 1, pp. 17--34, 2016.

\bibitem{Grenander}
\newblock U. Grenander,
\newblock \emph{A prediction problem in game theory},
\newblock Arkiv f\"or Matematik, vol. 3, pp. 371--379, 1957.


\bibitem{Hu}
\newblock W. Hu, and S. Dodelson,
\newblock \emph{Cosmic microwave background anisotropies},
\newblock Annual Review of Astronomy and Astrophysics, vol. 40, pp. 171--216, 2002.

\bibitem{Hurd}
\newblock H. L. Hurd, and A. Miamee,
\newblock \emph{Periodically Correlated random sequences: Spectral theory and practice},
\newblock Wiley Series in Probability and Statistics; Wiley Interscience. Hoboken, NJ: John Wiley and Sons, 2007.

\bibitem{Jones}
\newblock P. D. Jones,
\newblock \emph{Hemispheric surface air temperature variations: A reanalysis and an update to 1993},
\newblock Journal of Climate, vol. 7, pp. 1794-1802, 1994.

\bibitem{Kailath}
\newblock T. Kailath,
\newblock \emph{A view of three decades of linear filtering theory},
\newblock IEEE Transactions on Information Theory, Vol. 20, pp. 146--181, 1974.

\bibitem{Kakarala}
\newblock R. Kakarala,
\newblock \emph{The bispectrum as a source of phase-sensitive invariants for Fourier descriptors: A group-theoretic approach},
\newblock Journal of Mathematical Imaging and Vision, vol. 44, pp. 341--353, 2012.


\bibitem{Kallianpur}
\newblock G. Kallianpur, and V. Mandrekar,
\newblock \emph{Spectral theory of stationary H-valued processes},
\newblock J.  Multivariate Analysis, vol. 1, pp. 1--16, 1971.

\bibitem{Karhunen}
\newblock K. Karhunen,
\newblock \emph{Uber lineare Methoden in der Wahrscheinlichkeitsrechnung},
\newblock Annales Academiae Scientiarum Fennicae. Ser. A I, no. 37, 1947.

\bibitem{KassamPoor}
\newblock S. A. Kassam, and H. V. Poor,
\newblock \emph{Robust techniques for signal processing: A survey},
\newblock Proceedings of the IEEE, vol. 73, no. 3, pp. 433--481, 1985.

\bibitem{Kogo}
\newblock N. Kogo, and N.Komatsu,
\newblock \emph{Angular trispectrum of cmb temperature anisotropy from primordial non-Gaussianity with the full radiation
transfer function},
\newblock Phys. Rev., vol. D73, pp. 083007--083012, 2006.

\bibitem{Kolmogorov}
\newblock A. N. Kolmogorov,
\newblock \emph{Selected works by A. N. Kolmogorov. Vol. II: Probability theory and mathematical statistics. Ed. by A. N. Shiryayev},
\newblock Mathematics and its Applications. Soviet Series. 26. Dordrecht etc.
Kluwer Academic Publishers, 1992.


\bibitem{Luz2016}
\newblock  M. Luz, and M. Moklyachuk,
\newblock \emph{Estimation of Stochastic Processes with Stationary Increments and Cointegrated Sequences},
\newblock London:ISTE; Hoboken, NJ: John Wiley \& Sons, 282 p., 2019.



\bibitem{Marinucci}
\newblock D. Marinucci, and G. Peccati,
\newblock \emph{Random fields on the sphere},
\newblock London Mathematical Society Lecture Notes Series, vol. 389, Cambridge
University Press, Cambridge, 2011.


\bibitem {Moklyachuk:1981}
\newblock M. P. Moklyachuk,
\newblock \emph{Estimation of linear functionals of stationary stochastic processes and a two-person zero-sum game},
\newblock Stanford University Technical Report, no. 169, 87p., 1981.



\bibitem{Moklyachuk:2008}
\newblock M. P. Moklyachuk,
\newblock \emph{Robust estimations of  functionals of stochastic processes.},
\newblock Ky{\"\i}vsky\u\i\ Universytet, Ky{\"\i}v,\, 320p., 2008. (in Ukraine)

\bibitem{Moklyachuk:2008nonsm}
\newblock M. P. Moklyachuk,
\newblock \emph{Nonsmooth analysis and optimization},
\newblock Ky{\"\i}vsky\u\i\ Universytet, Ky{\"\i}v, \, 400p., 2008. (in Ukraine)

\bibitem{Moklyachuk:2015}
\newblock M. P. Moklyachuk,
\newblock \emph{Minimax-robust estimation problems for stationary stochastic sequences},
\newblock Statistics, Optimization \& Information Computing, vol. 3, no. 4, pp. 348--419, 2015.


\bibitem{Moklyachuk:2016}
\newblock M. Moklyachuk, and I. Golichenko,
\newblock \emph{Periodically correlated processes estimates},
\newblock LAP Lambert Academic Publishing, 2016.


\bibitem{Moklyachuk:2012}
\newblock M. Moklyachuk, and O. Masyutka,
\newblock \emph{Minimax-robust estimation technique for stationary stochastic processes},
\newblock LAP Lambert Academic Publishing, 2012.

\bibitem{MoklyachukYadrenko:1979}
\newblock M. P. Moklyachuk, and M. I. Yadrenko,
\newblock \emph{Linear statistical problems for homogeneous isotropic random fields on a sphere. I},
\newblock Theory of Probability and Mathematical Statistics, vol. 18, pp. 115-124, 1979.


\bibitem{MoklyachukYadrenko:1980}
\newblock M. P. Moklyachuk, and M. I. Yadrenko,
\newblock \emph{Linear statistical problems for homogeneous isotropic random fields on a sphere. II},
\newblock Theory of Probability and Mathematical Statistics, vol. 19, pp.129-139, 1980.


\bibitem{Muller}
\newblock C. M\"uller,
\newblock \emph{Spherical harmonics},
\newblock Lecture Notes in Mathematics 17. Berlin-Heidelberg-New York: Springer-Verlag, 1966.



\bibitem{Napolitano}
\newblock A. Napolitano,
\newblock \emph{Cyclostationarity: New trends and applications},
\newblock Signal Processing, vol. 120, pp. 385--408, 2016.


\bibitem{North}
\newblock G. R. North, and R. F. Cahalan,
\newblock \emph{Predictability in a solvable stochastic climate model},
\newblock J. Atmospheric Sciences, vol. 38, pp. 504-513, 1981.

\bibitem{Okamoto}
\newblock T. Okamoto, and W. Hu,
\newblock \emph{Angular trispectra of cmb temperature and polarization},
\newblock Phys. Rev., Vol. D66, p. 063008, 2002.



\bibitem{Rockafellar}
\newblock R. T. Rockafellar,
\newblock \emph{Convex Analysis},
\newblock Princeton University Press, 1997.


\bibitem{Rozanov}
\newblock Yu. A. Rozanov,
\newblock \emph{Stationary stochastic processes},
\newblock San Francisco-Cambridge-London-Amsterdam: Holden-Day, 1967.

\bibitem{Serpedin}
\newblock E. Serpedin, F. Panduru, I. Sari, and G. B. Giannakis,
\newblock  \emph{Bibliography on cyclostationarity},
\newblock Signal Processing, vol. 85, pp. 2233--2303, 2005.

\bibitem{SubbaRao2006}
\newblock T. Subba Rao and G. Terdik,
\newblock  \emph{Multivariate non-linear regression with applications},
\newblock In: P. Bertail, P. Doukhan, and P. Soulier (eds),
Dependence in Probability and Statistics. Springer Verlag, New York,
pp. 431-470, 2006.

\bibitem{SubbaRao2012}
\newblock T. Subba Rao and G. Terdik,
\newblock  \emph{Statistical analysis of spatio-temporal models and their applications},
\newblock In:  C. R. Rao (ed), Handbook of Statistics, Vol. 30, Elsevier B.V., pp. 521--541,  2012.


\bibitem{Terdik2015}
\newblock G. Terdik,
\newblock \emph{Angular spectra for non-Gaussian isotropic fields},
\newblock Brazilian Journal of Probability and Statistics,
vol. 29, no. 4, pp. 833--865, 2015.


\bibitem{Vastola}
\newblock  K. S. Vastola and H. V. Poor,
\newblock  \emph{ An analysis of the effects of spectral uncertainty on Wiener filtering},
\newblock  Automatica, vol. 28, pp. 289--293, 1983.


\bibitem{Wiener}
\newblock  N. Wiener,
\newblock \emph{Extrapolation, Interpolation and Smoothing of Stationary Time Series. With Engineering Applications},
\newblock  The M. I. T. Press, Massachusetts Institute of Technology, Cambridge, Mass., 1966.

\bibitem{Yadrenko}
\newblock M. I. Yadrenko,
\newblock \emph{Spectral theory of random fields},
\newblock Optimization Software Inc. Publications Division, New York, 1983.


\bibitem{Yaglom:1987a}
\newblock A. M. Yaglom,
\newblock \emph{Correlation theory of stationary and related random functions. Vol. 1: Basic results},
\newblock Springer Series in Statistics, Springer-Verlag, New York etc., 1987.

\bibitem{Yaglom:1987b}
\newblock A. M. Yaglom,
\newblock \emph{Correlation theory of stationary and related random functions. Vol. 2: Suplementary notes and references},
\newblock Springer Series in Statistics, Springer-Verlag, New York etc., 1987.



\end{thebibliography}
\end{document}